\documentclass[11pt]{article}
\usepackage[T1]{fontenc}
\usepackage{lmodern}
\usepackage[a4paper,margin=27mm]{geometry}
\usepackage{microtype}
\usepackage{amsmath,amssymb,amsthm,mathtools}
\usepackage[colorlinks=true,linkcolor=black,citecolor=black,urlcolor=blue]{hyperref}
\numberwithin{equation}{section}
\newtheorem{theorem}{Theorem}[section]
\newtheorem{proposition}[theorem]{Proposition}
\newtheorem{lemma}[theorem]{Lemma}
\newtheorem{corollary}[theorem]{Corollary}
\theoremstyle{definition}
\newtheorem{example}[theorem]{Example}
\theoremstyle{remark}
\newcommand{\R}{\mathbb R}

\newcommand{\HH}{\mathcal H}
\title{Sharp Bounds for the Nodal Measure of Clamped Plate \\Eigenfunctions in Smooth Domains
}
\author{Hairong Liu\textsuperscript{1}\qquad
        Long Tian\textsuperscript{2}\qquad
        Xiaoping Yang\textsuperscript{3}}
\hypersetup{
  pdftitle={Sharp Bounds for the Nodal Measure of Clamped Plate Eigenfunctions in Smooth Domains
    },
  pdfauthor={Hairong Liu, Long Tian, Xiaoping Yang}
}
\date{}
\begin{document}
\maketitle

\begin{abstract}
We investigate measure estimates for the nodal sets of
nonzero real eigenfunctions $u$ of the biharmonic operator, satisfying
$\Delta^2u=\lambda u$ in $\Omega$ with clamped boundary conditions
$u=\partial_\nu u=0$ on $\partial\Omega$. We establish the optimal
upper bound of order $\lambda^{1/4}$ for the $(n-1)$-dimensional
Hausdorff measure of their nodal sets on general bounded connected
$C^\infty$ domains.

Our approach combines an enhanced annular Carleman estimate 
with propagation of mass along a finite chain of
balls, these yield a doubling index bounded
 uniformly up to the boundary.
 Interior analyticity then gives the nodal measure away from the boundary. 
 To estimate the measure near the boundary,
we develop a subdivision argument using the boundary values of
$\Delta u$ and $\partial_\nu\Delta u$. Estimates from these data
reduce large doubling indices and control the
nodal measure in the region associated with a selected child cube
whenever both the doubling index of the parent cube and
$r\lambda^{1/4}$ are bounded, where $r$ is the side length of the
parent cube.
\end{abstract}

{\small
\par
}

\section{Introduction}

Let $\Omega\subset\R^n$ be a bounded connected $C^\infty$ domain.
We study the nodal sets of nonzero real eigenfunctions of the clamped
plate problem
\begin{equation}\label{main:problem}
  \begin{cases}
 \Delta^2u=\lambda u\quad\text{in }\Omega,\\[2mm]
 u=\partial_\nu u=0\quad\text{on }\partial\Omega,
  \end{cases}
\end{equation}
where $\nu$ is the outward unit normal. 

For Laplace eigenfunctions $-\Delta_g\phi=\lambda\phi$ on a compact
$n$-dimensional Riemannian manifold without boundary, Yau conjectured
the bounds
\begin{equation*}
 c\sqrt\lambda\le\mathcal H^{n-1}(\{\phi=0\})\le C\sqrt\lambda
\end{equation*}
for nonzero real eigenfunctions, with constants depending only on the
manifold. Donnelly and Fefferman \cite{DF88} proved this estimate in the
real-analytic case. Lin \cite{Lin91} developed an alternative approach
to upper bounds for general second-order elliptic equations in the
analytic setting.

On smooth surfaces, Donnelly and Fefferman \cite{DF90} and Dong
\cite{Dong} obtained the upper bound $C\lambda^{3/4}$, which Logunov
and Malinnikova \cite{LM} improved to $C\lambda^{3/4-\epsilon}$ for
some $\epsilon>0$. In arbitrary dimensions, Hardt and Simon \cite{HS}
established the upper bound $\lambda^{C\sqrt\lambda}$, and Logunov
\cite{LogUpper} later proved a polynomial upper bound. The sharp lower
bound in Yau's conjecture was proved by Logunov \cite{LogLower}; see
also the references therein for earlier results. More recently,
OpenAI \cite{OpenAI26} reported smooth counterexamples to the upper
bound on $S^3$ and $S^2\times\mathbb T^2$, with unbounded ratios of
nodal measure to $\sqrt\lambda$ along eigenfunction sequences for a
fixed metric.

For manifolds with boundary, Donnelly and Fefferman \cite{DFboundary}
established sharp upper and lower bounds in the real-analytic case.
Ariturk \cite{Ar2013} obtained lower bounds for Dirichlet and Neumann
eigenfunctions on smooth compact manifolds with boundary. For bounded
Euclidean domains, Logunov, Malinnikova, Nadirashvili, and Nazarov
\cite{LMNN} proved the sharp $C_\Omega\sqrt\lambda$ upper bound for
Dirichlet eigenfunctions on $C^1$ domains. Sharp Dirichlet bounds also
hold on quasiconvex Lipschitz domains \cite{ZZ} and bounded Lipschitz
polytopes \cite{CZ}, while the corresponding Neumann bound holds on
$C^{1,1}$ domains \cite{CY}.

For higher-order elliptic eigenvalue problems, Kukavica \cite{Kukavica}
proved the optimal nodal measure upper bound in the analytic case.
For the clamped plate,
Lin and Zhu \cite{LZ} proved the sharp $\lambda^{1/4}$ upper bound on
bounded real-analytic domains, together with sharp bounds for buckling,
Navier, and biharmonic Steklov problems. 
On smooth compact manifolds
without boundary, Zhu \cite{Zhu19} obtained doubling estimates for
$\Delta_g^2u=Wu$, $\|W\|_\infty\le M$, a polynomial bound for the mixed
zero set $\{u=\Delta_gu=0\}$ in dimensions $n\ge3$, and an implicit
bound for the full zero set. The subsequent work \cite{Zhu26} gives
the upper bound $CM^{\beta(n)}$, $\beta(n)>1/4$, for the full nodal
measure. For biharmonic Steklov problems on smooth compact manifolds
with boundary, Chang \cite{Chang} proved polynomial lower bounds for
interior and boundary nodal measure under regular-value assumptions.
For boundary doubling and spectral inequalities for clamped plates,
see \cite{MRV,LRR}. Tian and Yang \cite{TY14} used a frequency function
to obtain local nodal measure estimates for biharmonic functions.
In \cite{TY22}, they proved the sharp upper bound
\begin{equation*}
 \mathcal H^{n-1}(\{u=0\}\cap\Omega)\le C\lambda^{1/4}
\end{equation*}
for eigenfunctions satisfying $\Delta^2u=\lambda u$ on bounded
$C^{(n+1)/2}$ domains whose boundaries are piecewise analytic and
analytic outside a finite union of compact $(n-2)$-dimensional
submanifolds, under the boundary-operator hypotheses in \cite{TY22}.

We prove the following sharp bound on smooth domains.
\begin{theorem}\label{main:sharp}
Let $\Omega\subset\mathbb{R}^n$ ($n\geq2$) be a bounded connected domain
with $C^{\infty}$ boundary. Let $u\in H^4(\Omega)\cap H^2_0(\Omega)$
be a nonzero real eigenfunction of \eqref{main:problem}. Then
there is $C_\Omega>0$ such that 
\begin{equation}\label{main:bound}
 \HH^{n-1}(\{u=0\}\cap\Omega)\le C_\Omega\lambda^{1/4}.
\end{equation}
\end{theorem}

\vspace{3mm}

Theorem~\ref{main:sharp} gives the sharp upper bound
$C_\Omega\lambda^{1/4}$ for the $(n-1)$-dimensional Hausdorff measure
of clamped plate nodal sets on every bounded connected $C^\infty$
domain $\Omega\subset\R^n$, $n\ge2$. Unlike the real-analytic case,
the eigenfunctions need not extend analytically across a smooth
boundary, although they are analytic in $\Omega$. The main difficulty
is to control the nodal measure arbitrarily close to $\partial\Omega$.

Set $k=\lambda^{1/4}$. An annular Carleman estimate controlling the
conjugates of $-\Delta\pm k^2$, combined with a singular estimate
from clamped reflection and coefficient freezing, gives local
doubling uniformly up to $\partial\Omega$. The estimates are matched
at a scale where coefficient errors can be absorbed. Propagation of
$L^2$ mass along a finite ball chain bounds all local doubling
indices by $C_\Omega(1+k)$; holomorphic extension then yields the
interior nodal measure estimates.

Near the boundary, we subdivide tangential cubes $Q$ and set
$\mu=k\ell(Q)$. Propagation from the traces $\Delta u$ and
$\partial_\nu\Delta u$ selects a child with a halved doubling index
when the parent index exceeds $K(1+\mu)$, and controlled nodal measure
when both the parent index and $\mu$ are bounded. Together with the
decay of $\mu$, these alternatives give a contracting weighted sum
and hence the boundary estimate.

Section~\ref{sec2} contains the regularity and propagation estimates.
Sections~\ref{sec3} and~\ref{sec4} prove the Carleman and doubling
estimates. Section~\ref{sec5} proves Theorem~\ref{main:sharp}
and its sharpness on the unit ball.

\section{Regularity and propagation}\label{sec2}

After translation and rotation at $p\in\partial\Omega$, write locally
\begin{equation*}
  \partial\Omega=\{x_n=g(x')\},\qquad
  \Omega=\{x_n>g(x')\},
\end{equation*}
where $g$ is smooth, $g(0)=0$, and $\nabla g(0)=0$. Set
\begin{equation*}
  F(y,t)=(y,g(y)+t),\qquad \det DF=1.
\end{equation*}
Choose $a_0>0$ so that $F$ is a boundary chart on $(-a_0,a_0)^n$.
For $0<a<a_0$, write
\begin{equation*}
  Q_a=(-a,a)^n,\qquad
  Q_a^+=(-a,a)^{n-1}\times(0,a),\qquad
  D_a=F(Q_a^+).
\end{equation*}
For $U=u\circ F$, the clamped conditions on the graph are
equivalent to
\begin{equation*}
  U(y,0)=\partial_tU(y,0)=0.
\end{equation*}
Fix a finite boundary atlas with uniform bounds on the maps and their
inverses at each fixed derivative order. Boundary conditions below
refer only to the physical boundary.

For fixed $0<a<b<a_0$ and an integer $j\geq0$, the clamped elliptic
estimate of Agmon, Douglis, and Nirenberg \cite{ADN} gives
\begin{equation}\label{nod:boundary-regularity}
 \|u\circ F\|_{H^{j+4}(Q_a^+)}
 \le C\bigl(
   \|(\Delta^2u)\circ F\|_{H^j(Q_b^+)}
   +\|u\circ F\|_{L^2(Q_b^+)}\bigr).
\end{equation}
Here $C$ depends on $n,j,a,b$ and finitely many seminorms of $g$,
uniformly over the chosen charts.

For $\Delta^2u=\lambda u$, iterating
\eqref{nod:boundary-regularity} and applying Sobolev embedding give,
for each integer $q\geq0$,
\begin{equation}\label{nod:finite-regularity}
  \|u\circ F\|_{C^q(Q_a^+)}
  \leq C_q(1+k)^{b_q}
       \|u\circ F\|_{L^2(Q_b^+)},
\end{equation}
where $k=\lambda^{1/4}$ and $C_q,b_q>0$ are independent of
$\lambda$. In charts rescaled by $x=x_0+ry$, the same estimate holds
with $k$ replaced by $kr$.

For the interior estimates, set
\begin{equation}\label{nod:lift-operator}
 P=\Delta_x^2+\partial_s^4.
\end{equation}

For $\Delta^2u=k^4u$, the lift
$h(x,s)=e^{e^{i\pi/4}ks}u(x)$ satisfies $Ph=0$.
For $Ph=0$ on a ball $ B\subset\R^{n+1}$ and a concentric
ball $ B'\Subset B$, analytic elliptic regularity
\cite{MN} gives
\begin{equation}\label{nod:analytic-regularity}
 \sup_{ B'}|\partial^\beta h|
 \le C_0 C_1^{|\beta|}|\beta|!\,\|h\|_{L^2( B)}.
\end{equation}
Hence $h$ extends holomorphically to a fixed complex neighborhood
of each smaller ball, with bound $C\|h\|_{L^2(B)}$.
Applying \cite[Theorem~4]{AEWZ2014} to $\operatorname{Re}h$ and
$\operatorname{Im}h$ along a finite ball chain gives the following lemma.

\begin{lemma}
\label{nod:analytic-propagation}
Let $ B\subset\R^{n+1}$ be a ball, let
$ E\Subset B$ be a nonempty ball, and let
$K\Subset B$ be compact. For $Ph=0$ in $ B$,
\begin{equation}\label{nod:analytic-propagation-estimate}
 \|h\|_{L^2(K)}
 \le C\|h\|_{L^2( B)}^{1-\alpha}
       \|h\|_{L^2( E)}^\alpha
\end{equation}
with constants $C>0$ and $\alpha\in(0,1)$ depending only on $n,B,E,K$.
\end{lemma}

For the boundary estimate, let $g_p$ be the graph at $p\in\partial\Omega$
and set $g(y')=\ell^{-1}g_p(\ell y')$. For sufficiently small $\ell$,
\begin{equation*}
 g\in C^4(B_4^{n-1}),\qquad g(0)=0,\qquad
 \|\nabla g\|_{L^\infty(B_4^{n-1})}\le\frac1{20},\qquad
 \|g\|_{C^4(B_4^{n-1})}\le G,
\end{equation*}
with $G$ independent of $p$ and $\ell$. Put
\begin{equation*}
 D=B_4^n\cap\{x_n>g(x')\},\qquad
 \Gamma=B_3^n\cap\{x_n=g(x')\},\qquad
 E_{\rm in}=B_{1/16}^n(\tfrac12e_n).
\end{equation*}

\begin{lemma}\label{lem-cauchy}
Let $u\in H^4(D)$ satisfy $\Delta^2u=\lambda u$ in $D$,
with $\lambda\geq0$, and
\begin{equation}\label{nod:clamped-local}
  u=\partial_\nu u=0
  \qquad\text{on }\partial D\cap B_4^n.
\end{equation}
Set
\begin{equation*}
  T=\|\Delta u\|_{L^2(\Gamma)}
    +\|\partial_\nu\Delta u\|_{L^2(\Gamma)}.
\end{equation*}
There exist constants $C>0$ and $\alpha\in(0,1)$,
depending only on $n$ and $G$, such that, for every
$M\geq\|u\|_{L^2(D)}$,
\begin{equation}\label{nod:cauchy-estimate}
  \|u\|_{L^2(E_{\mathrm{in}})}
  \leq C e^{C\lambda^{1/4}} M^{1-\alpha}T^\alpha.
\end{equation}
The same estimate holds with $E_{\mathrm{in}}$ replaced
by any ball $E\Subset D\cap B_1^n$, with constants
depending additionally on positive lower bounds for
the radius of $E$ and its distance from $\partial D$.
\end{lemma}
\begin{proof}
The case $T\geq M$ is immediate, so assume $T<M$.
Write $B_r=B_r^{n+1}(0)$ and set
\begin{equation*}
  U(x,s)=e^{e^{i\pi/4}\lambda^{1/4}s}u(x),
\end{equation*}
so that $PU=0$. Let $\widetilde U$ be its zero extension
across the cylindrical graph $\{x_n=g(x')\}\times\R_s$.
The clamped conditions \eqref{nod:clamped-local}
imply that $\widetilde U\in H^2(B_2)$.

For $\varphi\in C_c^\infty(B_2)$, Green's formula gives
\begin{equation*}
  \langle P\widetilde U,\varphi\rangle
  =
  \int_{\Gamma\times\mathbb R}
  e^{e^{i\pi/4}\lambda^{1/4}s}
  \bigl(
    \Delta u\,\partial_\nu\varphi
    -\partial_\nu\Delta u\,\varphi
  \bigr)\,dS_x\,ds.
\end{equation*}
The trace theorem gives
\begin{equation*}
  \|\widetilde U\|_{L^2(B_2)}
  \leq C e^{C\lambda^{1/4}}M,
  \qquad
  \|P\widetilde U\|_{(H_0^2(B_2))^*}
  \leq C e^{C\lambda^{1/4}}T.
\end{equation*}
Coercivity on $H_0^2(B_2)$ gives $w\in H_0^2(B_2)$ with
\begin{equation*}
  Pw=P\widetilde U\quad\text{in }B_2,
  \qquad
  \|w\|_{H^2(B_2)}
  \leq C e^{C\lambda^{1/4}}T.
\end{equation*}
Thus $h=\widetilde U-w$ satisfies $Ph=0$ in $B_2$.
The slope bound on $g$ ensures that
\begin{equation*}
  E_{\mathrm{out}}
  =B_{1/8}^{n+1}\bigl((-\tfrac12 e_n,0)\bigr)
\end{equation*}
lies below the graph, where $h=-w$. Since $T<M$,
\begin{equation*}
  \|h\|_{L^2(E_{\mathrm{out}})}
  \leq C e^{C\lambda^{1/4}}T,
  \qquad
  \|h\|_{L^2(B_2)}
  \leq C e^{C\lambda^{1/4}}M.
\end{equation*}
Applying Lemma~\ref{nod:analytic-propagation} with
$E=E_{\mathrm{out}}$ and
$K=\overline{B_{3/2}}$ gives
\begin{equation*}
  \|h\|_{L^2(B_{3/2})}
  \leq C e^{C\lambda^{1/4}}M^{1-\alpha}T^\alpha.
\end{equation*}
Since $\widetilde U=h+w$ and $T\leq M^{1-\alpha}T^\alpha$,
\begin{equation*}
  \|\widetilde U\|_{L^2(B_{3/2})}
  \leq C e^{C\lambda^{1/4}}M^{1-\alpha}T^\alpha.
\end{equation*}
Since $E_{\mathrm{in}}\times(0,\tfrac12)\subset B_{3/2}$
and $|U(x,s)|\geq|u(x)|$ for $s\geq0$,
\begin{equation*}
  \|u\|_{L^2(E_{\mathrm{in}})}
  \leq \sqrt{2}\,
       \|\widetilde U\|_{L^2(B_{3/2})}
  \leq C e^{C\lambda^{1/4}}M^{1-\alpha}T^\alpha.
\end{equation*}
For the other target balls, use
$E\times(0,\tfrac12)\subset B_{3/2}$ in the last step.
\end{proof}

At scale $\ell$, the exponential factor in Lemma~\ref{lem-cauchy}
is $e^{Ck\ell}$, with $C$ and $\alpha$ uniform in $p$ and small
$\ell$ for fixed rescaled observation sets.

\section{The annular Carleman estimate}\label{sec3}

For $x_0\in\overline\Omega$ and $0<R\le R_0$, write
\begin{equation*}
 \Omega_{x_0,R}=\{y:x_0+Ry\in\Omega\},\qquad r=|y|,
 \qquad T(r)=\log\left(\frac{r_c}{r}\right),\qquad r_c=3e^{T_0}.
\end{equation*}
Choose $a>0$, then $T_0\gg1$ and $0<b\ll a$.
All enlarged supports lie in $\{r<4\}$, where $T>0$.
Let $d$ be the signed distance to $\partial\Omega$, positive in $\Omega$.

If $\operatorname{dist}(x_0,\partial\Omega)\le4R$, set
\begin{equation}\label{car:eq:phase} d_R(y)=R^{-1}d(x_0+Ry),\qquad
 \Phi(y)=T(r)-a\log T(r)+
              b\frac{d_R(y)-d_R(0)}{rT(r)^2}.
\end{equation}
If $\operatorname{dist}(x_0,\partial\Omega)>4R$, we
set $\Phi(y)=T(r)-a\log T(r)$.  

\begin{theorem}
\label{car:thm:annular}
Let $\Omega\subset\mathbb R^n$, $n\geq2$, be a bounded
$C^\infty$ domain. There exist an integer $M\geq1$ and constants
$R_0>0$ and $C\geq1$, all depending only on $n$ and $\Omega$,
for which the following estimate holds.

For every
\begin{equation*}
  x_0\in\overline{\Omega},\qquad
  0<R\leq R_0,\qquad
  0<\delta<1,\qquad
  k\geq0,
\end{equation*}
set
\begin{equation*}
  T_*=\log\left(\frac{r_c}{\delta}\right),\qquad r_c=3e^{T_0}.
\end{equation*}
Suppose that $w\in C^\infty(\overline{\Omega_{x_0,R}})$ is supported in
$\overline{\Omega_{x_0,R}}\cap\{\delta<|y|<3\}$ and satisfies
\begin{equation*}
  w=\partial_{\nu} w=0
  \qquad\text{on }\partial\Omega_{x_0,R}.
\end{equation*}
Then
\begin{equation}\label{car:eq:annular-scalar}
  \tau^3
  \|e^{\tau\Phi}w\|_{L^2(\Omega_{x_0,R})}
  \leq C T_*^M
  \bigl\|r^4e^{\tau\Phi}(\Delta^2-k^4)w
  \bigr\|_{L^2(\Omega_{x_0,R})},
\end{equation}
provided
\begin{equation*}
  \tau\geq C\big(1+k+T_*^M\big).
\end{equation*}
\end{theorem}

\vspace{3mm}

For the radial part $F(T)=T-a\log T$,
\begin{equation*}
 F'=1-a/T\in[3/4,1],\qquad F''=a/T^2,\qquad
 |\Phi-F(T)|\le bT^{-2}.
\end{equation*}
The oscillation on each fixed-ratio annulus is bounded, and small $b$ gives
\begin{equation*}
 \inf_{0<r\le1}\Phi-\sup_{2<r<3}\Phi\ge g>0.
\end{equation*}

\begin{lemma}
\label{car:lem:phase}
On a fixed normalized annulus $y=\rho z$, $0<\rho\leq4$,
set
\begin{equation*}
  t=\log\left(\frac{r_c}{\rho}\right),\qquad \eta=t^{-2},\qquad
  \phi(z)=\Phi(\rho z),\qquad h=\tau^{-1},\qquad \theta=hk\rho.
\end{equation*}
All required positive-order derivatives of $\phi$
are uniformly bounded.

There exist $c_*,c,C>0$ such that, if
$\theta^2\leq c_*\eta$, the symbols
\begin{equation*}
  q_{\pm}
  =|\xi|^2-|\nabla\phi|^2\pm\theta^2
   +2i\xi\cdot\nabla\phi,
\end{equation*}
satisfy, on a fixed bounded frequency box,
\begin{equation}\label{car:eq:second-completion}
  C\eta^{-1}(\operatorname{Re}q_{\pm})^2
  +\{\operatorname{Re}q_{\pm},\operatorname{Im}q_{\pm}\}
  \geq c\eta.
\end{equation}
Here $\{A,B\}=\partial_\xi A\cdot\partial_z B
-\partial_z A\cdot\partial_\xi B$ denotes the Poisson bracket
in the variables $(z,\xi)$.

If $R\rho\leq\varepsilon_b\eta$, with
$\varepsilon_b>0$ sufficiently small, then
\begin{equation}\label{car:eq:normal-margin}
  \partial_n\phi\geq c_b\eta
\end{equation}
on the boundary, where $\partial_n$ points inward.
After decreasing $c_b$, the bound holds in a boundary collar of width
$c_b'\eta$ in the normalized $z$-coordinates; all constants are
independent of $x_0$ and $\rho$.
\end{lemma}

\begin{proof}
Write
\begin{equation*}
  s=|z|,\qquad T=t-\log s,\qquad p=1-a/T.
\end{equation*}
On the fixed normalized annulus, $T\sim t$.
The signed-distance bounds give, for every fixed required order $j$,
\begin{equation*}
  \|\phi-F(T)\|_{C^j}\leq C_jb\eta.
  \end{equation*}

At a boundary point of radius $r$, put $d_0=d_R(0)$.
Then
\begin{equation*}
  0\leq d_0\leq r,\qquad
  \partial_n r=-\frac{d_0}{r}+O(Rr).
\end{equation*}
Direct differentiation gives
\begin{equation*}
  \partial_n\Phi
  =
  p\frac{d_0}{r^2}
  +\frac{b}{rT^2}
   \left[1-\left(1-\frac2T\right)\frac{d_0^2}{r^2}\right]
  +O(R).
\end{equation*}
Since $p\ge3/4$ and $0\le d_0/r\le1$,
\begin{equation*}
  \partial_n\Phi
  \geq c\left(\frac{d_0}{r^2}+\frac{b}{rT^2}\right)-CR.
\end{equation*}
Thus $r\sim\rho$, $T^{-2}\sim\eta$, and
$R\rho\le\varepsilon_b\eta$ imply \eqref{car:eq:normal-margin}.
The uniform Hessian bound extends it to this collar.
Rescaling fixed cutoffs in the normal variable gives cutoffs $\chi$
supported in this collar, with
\begin{equation*}
 \|\partial_z^\alpha\chi\|_{L^\infty}
 \le C_\alpha\eta^{-|\alpha|}.
\end{equation*}

The same differentiation gives
\begin{align*}
  \frac{1}{4}
  \{\operatorname{Re}q_{\pm},\operatorname{Im}q_{\pm}\}
  &=
  \frac{ap^2}{s^4T^2}
  +\frac{2p+a/T^2}{s^4}(\xi\cdot z)^2
  -\frac p{s^2}\operatorname{Re}q_{\pm}
  {\pm}\frac{ p}{s^2}\theta^2
  +O(b\eta)\\
  &\geq ca\eta-Cb\eta-C\theta^2
          -C|\operatorname{Re}q_{\pm}|.
\end{align*}
For $b\ll a$ and $\theta^2\le c_*\eta$, Young's inequality gives
\begin{equation*}
  \{\operatorname{Re}q_{\pm},\operatorname{Im}q_{\pm}\}
  \geq c\eta-C|\operatorname{Re}q_{\pm}|
  \geq \frac c2\eta
       -C'\eta^{-1}(\operatorname{Re}q_{\pm})^2.
\end{equation*}
The Poisson bracket is invariant under the change to boundary normal
coordinates.
\end{proof}

\begin{proof}[Proof of Theorem~\ref{car:thm:annular}]
Put $D=-i\partial$ and define
\begin{equation}\label{car:eq:global-factors}
 Q_{\pm,\tau}=e^{\tau\Phi}(-\Delta\pm k^2)e^{-\tau\Phi},
 \qquad
 P_\tau=Q_{-,\tau}Q_{+,\tau}
       =e^{\tau\Phi}(\Delta^2-k^4)e^{-\tau\Phi}.
\end{equation}
For $v$ supported in $\{\delta<r<3\}$, define
\begin{equation}\label{car:eq:global-energy}
 \mathcal E_\tau(v)^2=
 \sum_{j=0}^4\tau^{6-2j}\|r^jD^jv\|_2^2
 +\sum_{\epsilon=\pm}\sum_{j=0}^2
       \tau^{3-2j}\|r^{2+j}D^jQ_{\epsilon,\tau}v\|_2^2.
\end{equation}
Here $D^j$ denotes all derivatives of order $j$, and the norms are
taken over $\Omega_{x_0,R}$.
For $v=e^{\tau\Phi}w$, we will prove the stronger bound
\begin{equation}\label{car:eq:annular-enhanced}
 \mathcal E_\tau(v)
 \le C T_*^M\|r^4P_\tau v\|_{L^2(\Omega_{x_0,R})}
\end{equation}
under the stated lower bound on $\tau$.

\par\medskip\noindent\textbf{Step 1. Tangential calculus and boundary coordinates.}\par

For $\phi(z)=F(T)$,
\begin{equation}\label{car:eq:spectral-restriction}
 \left.\{\operatorname{Re}q_-,\operatorname{Im}q_-\}\right|_{q_-(z,\xi)=0}
 =4\left(\frac{ap^2}{s^4T^2}-\frac{p\theta^2}{s^2}\right).
\end{equation}
Since
\begin{equation*}
 \sup_{0<\rho\le4}\rho\log(r_c/\rho)<\infty,
 \qquad \tau\ge Ck\ \Longrightarrow\ \theta^2\le c_*\eta,
\end{equation*}
the spectral restriction is uniform in $\rho$.

On $X=\mathbb R^{n-1}\times\mathbb R_+$, fix
\begin{equation*}
 \vartheta\in C_c^\infty(\mathbb R^{n-1};\mathbb R),
 \qquad \vartheta(z)=\vartheta(-z),
 \qquad \vartheta=1\quad\text{near }0,
\end{equation*}
and define the tangential Weyl quantization by
\begin{equation*}
 \begin{aligned}
  (\operatorname{Op}_h(a)v)(x',x_n)
  ={}&\frac{1}{(2\pi h)^{n-1}}
  \iint_{\mathbb R^{2(n-1)}}
  e^{i(x'-y')\cdot\xi'/h}\vartheta(x'-y')\\
  &\qquad\times a\!\left(\frac{x'+y'}2,x_n,\xi'\right)
  v(y',x_n)\,dy'\,d\xi'.
 \end{aligned}
\end{equation*}
Symbols are smoothly extended beyond the relevant supports.
For $a\in S^m_{1,0}$, the bounds used below are
\begin{equation*}
 |\partial_x^\alpha\partial_{\xi'}^\beta a(x,\xi')|
 \le C_{\alpha\beta}\eta^{-M_{\alpha\beta}}
       \langle\xi'\rangle^{m-|\beta|},
 \qquad\text{for all }\alpha,\beta.
\end{equation*}
All constants and exponents in this calculus are independent of $h$ and $\eta$.

For integers $j\ge0$, we use the semiclassical Sobolev norms
\cite[Section~7.1, Eq.~(7.1.1)]{Zworski},
\begin{equation*}
 \|v\|_{H_h^j}^2=\sum_{|\alpha|\le j}\|(hD)^\alpha v\|_2^2.
\end{equation*}
For $s\in\mathbb R$, define the tangential norms by
\begin{equation*}
 \|v\|_{H_h^s(\mathrm{tan})}^2
 :=\int_0^\infty
 \|\langle hD'\rangle^s v(\cdot,x_n)\|_{L^2(\mathbb R^{n-1})}^2\,dx_n,
 \qquad D'=-i\nabla_{x'}.
\end{equation*}
For $a\in S^{m_1}_{1,0}$ and $b\in S^{m_2}_{1,0}$, the tangential
calculus \cite[Theorems~4.23, 9.5, and~9.6]{Zworski} gives
\begin{equation*}
 \|\operatorname{Op}_h(a)v\|_{H_h^s(\mathrm{tan})}
 \le C_s\eta^{-K_s}\|v\|_{H_h^{s+m_1}(\mathrm{tan})}.
\end{equation*}
For the Weyl product $a\#b$ and every integer $N\ge1$,
\begin{equation*}
 \begin{aligned}
  \operatorname{Op}_h(a)\operatorname{Op}_h(b)
   &=\operatorname{Op}_h(a\#b)+\mathcal R_h,\\
  a\#b&=\sum_{j=0}^{N-1}h^jc_j+h^Nr_N,\qquad c_0=ab,\\
  c_j&\in S^{m_1+m_2-j}_{1,0},\qquad
  r_N\in S^{m_1+m_2-N}_{1,0}.
 \end{aligned}
\end{equation*}
The seminorms of $c_j$ and $r_N$ satisfy analogous polynomial bounds
in $\eta^{-1}$, uniformly in $x_n$. The remainder $\mathcal R_h$
introduced by the kernel cutoff $\vartheta(x'-y')$ satisfies,
for every integer $L\ge0$ and $q\ge0$,
\begin{equation*}
 \|\mathcal R_h\|_{H_h^{-q}(\mathrm{tan})\to H_h^q(\mathrm{tan})}
 \le C_{L,q}h^L\eta^{-K_{L,q}}.
\end{equation*}
Normal derivatives are counted separately.
For a real scalar symbol $a\in S^m_{1,0}$ with $a\ge0$, the sharp
G\aa rding inequality \cite[Theorem~9.11]{Zworski} gives
\begin{equation}\label{car:eq:garding}
 \operatorname{Re}(\operatorname{Op}_h(a)w,w)
 \ge-Ch\eta^{-K}\|w\|_{H_h^{(m-1)/2}(\mathrm{tan})}^2.
\end{equation}
For a pointwise nonnegative Hermitian matrix symbol
$a\in S^0_{1,0}$ we use
\begin{equation*}
 \operatorname{Re}(\operatorname{Op}_h(a)W,W)\ge-Ch\eta^{-K}\|W\|_2^2.
\end{equation*}
The matrix estimate follows from positive Gaussian quantization and its $O(h)$
comparison with Weyl quantization, uniformly up to $x_n=0$.

Let
\begin{equation}\label{car:eq:word-class}
 \mathcal D=\sum_{a=0}^r\operatorname{Op}_h(t_a)(hD_n)^a,
 \qquad t_a\in S^{r-a}_{1,0},\qquad 0\le r\le4.
\end{equation}
Let $\chi,\psi\in S^0_{1,0}$ be phase-space cutoff symbols satisfying
these seminorm bounds and
\begin{equation*}
 \operatorname{dist}\bigl(\operatorname{supp}\chi,
             \operatorname{supp}(1-\psi)\bigr)\ge c\eta^M,
\end{equation*}
uniformly in $x_n\ge0$, where supports and distance are taken in
the tangential phase space $(x',\xi')$, and $c>0$, $M>0$ are
independent of $h$ and $\eta$.
For each fixed $N$, and every integer $s\ge0$ with $r+s\le4$,
\begin{equation}\label{car:eq:exterior}
 \|(I-\operatorname{Op}_h\psi)\mathcal D\operatorname{Op}_h\chi\,v
          \|_{H_h^s(X)}
 \le C_Nh^N\eta^{-K_N}\|v\|_{H_h^4(X)}.
\end{equation}
We use the identities
\begin{equation*}
 [hD_n,\operatorname{Op}_h(a)]
   =-ih\operatorname{Op}_h(\partial_na),\qquad
 \partial^\alpha(1-\psi)\,\partial^\beta\chi=0.
\end{equation*}
The composition formula and tangential disjoint-support estimates
\cite[Theorems~4.25 and~9.5--9.7]{Zworski} yield
\eqref{car:eq:exterior}, with at most $r+s\le4$ normal derivatives
of $v$ and constants uniform in $x_n\ge0$.

Let $\mathcal L$ denote $-\Delta_z$ expressed in boundary normal
coordinates $x=(x',x_n)$, where $x_n$ is the inward distance to the
boundary in the normalized $z$-coordinates. Then
\begin{equation*}
 \mathcal L=-\partial_n^2-b_0(x)\partial_n+B(x_n,x',D'),
\end{equation*}
where $g_t$ is the tangential metric, $t=x_n$, and
$J=\sqrt{\det g_t}$. Here
\begin{equation*}
 b_0=\partial_t\log J,\qquad
 B=-J^{-1}\sum_{\alpha,\beta<n}
       \partial_\alpha(Jg_t^{\alpha\beta}\partial_\beta).
\end{equation*}
With $m_0=(J(x',t)/J(x',0))^{-1/2}$,
$2\partial_t m_0+b_0m_0=0$, so
\begin{equation}\label{car:eq:gauge}
 m_0^{-1}\mathcal Lm_0=Q=D_n^2+A,\qquad
 A=m_0^{-1}Bm_0+\tfrac12\partial_t b_0+\tfrac14b_0^2.
\end{equation}
The principal symbol $a_2(x,\xi')=g_t^{\alpha\beta}\xi_\alpha\xi_\beta$
of $A$ is real positive definite.

For the shape operator $S=-d\nu_{\rm in}$, rescaling the boundary
graph gives, at each fixed required order,
\begin{equation*}
 |\partial^\alpha S|\le C_\alpha R\rho.
\end{equation*}
The identities
\begin{equation*}
 g_t=g_0\bigl((I-tS)\,\cdot,(I-tS)\,\cdot\bigr),
 \qquad \frac{J(x',t)}{J(x',0)}=\det(I-tS)
\end{equation*}
hold in a fixed boundary normal coordinate neighborhood of the
normalized domain, with $I-tS$ uniformly invertible.
Writing $A=\sum_{|\beta|\le2}a_\beta(x)D_{x'}^\beta$, we obtain
\begin{equation}\label{car:eq:curvature-scale}
 |\partial_x^\alpha\partial_n a_\beta(x)|
 \le C_{\alpha\beta}\kappa,\qquad \kappa\le CR\rho,
\end{equation}
uniformly in this neighborhood, for each fixed $\alpha$.
The coordinate maps and $m_0^{\pm1}$ have uniform bounds at each
fixed derivative order; multiplication by $m_0^{\pm1}$ preserves
the clamped boundary conditions.

Set $\mu=k\rho$, $\theta=h\mu$, and
\begin{align}\label{car:eq:local-factors}
 D_h&=hD_n+i\phi_n,
 &A_h&=e^{\phi/h}h^2Ae^{-\phi/h},\notag\\
 Q_{j,h}&=D_h^2+A_h+(-1)^j\theta^2,
 &P_h&=Q_{1,h}Q_{2,h},\qquad j=1,2.
\end{align}
Then
\begin{equation*}
 [Q_{1,h},Q_{2,h}]=0,\qquad Q_{1,h}-Q_{2,h}=-2\theta^2.
\end{equation*}
In the boundary collar $0\le x_n<c_b'\eta$ of
Lemma~\ref{car:lem:phase}, set
\begin{equation*}
 f=\phi_n\ge c\eta,\qquad
 a=a_2(x,\xi'+i\nabla'\phi),\qquad z_j=a+(-1)^j\theta^2.
\end{equation*}

\par\medskip\noindent\textbf{Step 2. Normal roots and factorization.}\par

Write $s_j=\sqrt{z_j}=p_j+i\ell_j$, with $p_j\ge0$, for $j=1,2$.
Fix a sufficiently large $L$. For each fixed $x_n$, we localize in
the tangential phase space $(x',\xi')$. On $|\xi'|\le2L$, each factor
has two overlapping regions:
\begin{equation*}
 \begin{array}{ll}
 \text{(i)\ regular:}&p_j>f/4\ge c\eta,\\[2mm]
 \text{(ii)\ stable:}&p_j<f/2,\quad
          \operatorname{Im}(-if\pm is_j)=-f\pm p_j\le-f/2.
 \end{array}
\end{equation*}
For $p_j>f/4$, define the first-order scalar symbols
\begin{equation*}
 \mathfrak a_j=\xi_n+if+is_j,\qquad
 \mathfrak b_j=\xi_n+if-is_j=H+iK,\qquad
 H=\xi_n+\ell_j,\quad K=f-p_j.
\end{equation*}
The root of $\mathfrak a_j$ satisfies
$\operatorname{Im}(-if-is_j)=-(f+p_j)\le-c\eta$.
For $\mathfrak b_j$, we claim
\begin{equation}\label{car:eq:first-completion}
 C\eta^{-9}K^2+\{H,K\}\ge c\eta.
\end{equation}
To construct the partition without differentiating square roots
near zero, set $\gamma=3/8$ and
\begin{equation}\label{car:eq:branch-partition}
 \begin{aligned}
  G_j&=(\operatorname{Im}z_j)^2
       +4\gamma^2f^2\operatorname{Re}z_j-4\gamma^4f^4\\[2mm]
     &=4(p_j^2-\gamma^2f^2)(\ell_j^2+\gamma^2f^2).
 \end{aligned}
\end{equation}
Since
\begin{equation*}
 G_j|_{p_j=f/4}\le-c\eta^4,\qquad G_j|_{p_j=f/2}\ge c\eta^4,
\end{equation*}
applying fixed smooth cutoff functions to $G_j/(c'\eta^4)$ gives
an overlapping partition with seminorms polynomial in $\eta^{-1}$.
Taking products of the two partitions gives cutoff symbols supported
in the four intersections
\begin{equation*}
 \begin{array}{ll}
  \{p_1>f/4\}\cap\{p_2>f/4\},&
  \{p_1>f/4\}\cap\{p_2<f/2\},\\[2mm]
  \{p_1<f/2\}\cap\{p_2>f/4\},&
  \{p_1<f/2\}\cap\{p_2<f/2\}.
 \end{array}
\end{equation*}

On $\{p_j>f/4\}$, the chain rule gives, for $|\alpha|\ge1$,
\begin{equation}\label{car:eq:root-seminorms}
 |\partial^\alpha s_j|\le C_\alpha\eta^{-(2|\alpha|-1)},
 \qquad
 |\partial^\alpha(s_j^{-1})|
       \le C_\alpha\eta^{-(2|\alpha|+1)}.
\end{equation}
Write $\mathfrak a_j=\alpha_j+i\beta_j$, where
$\alpha_j=\xi_n-\ell_j$ and $\beta_j=f+p_j$ are real-valued.
On $H=0$, $\alpha_j=-2\ell_j$, $c\eta\le|\mathfrak a_j|\le C$, and
\begin{align}\label{car:eq:product-bracket}
 \{\operatorname{Re}(\mathfrak a_j\mathfrak b_j),
   \operatorname{Im}(\mathfrak a_j\mathfrak b_j)\}
 &=|\mathfrak a_j|^2\{H,K\}
 +K\bigl(\alpha_j\{H,\alpha_j\}+\beta_j\{H,\beta_j\}\nonumber\\[2mm]
 &+\alpha_j\{K,\beta_j\}-\beta_j\{K,\alpha_j\}\bigr)
       +K^2\{\alpha_j,\beta_j\}.
\end{align}
Differentiate before restricting to $H=0$. Since
$\operatorname{Re}(\mathfrak a_j\mathfrak b_j)=-\beta_jK$, \eqref{car:eq:second-completion} and
\eqref{car:eq:root-seminorms} give
\begin{equation*}
 \{H,K\}\ge c\eta-C\eta^{-4}|K|-C\eta^{-4}K^2.
\end{equation*}
Young's inequality gives \eqref{car:eq:first-completion} on $H=0$.
As $\{H,K\}$ is independent of $\xi_n$, it holds throughout the
set $\{p_j>f/4\}$; choose the fixed frequency box to contain
$\xi_n=-\ell_j$.

For each $j=1,2$, on $\{p_j>f/4\}$ choose
$S_j=\operatorname{Op}_h(s_j+hs_{j,0})$, where the order-zero
correction $s_{j,0}$ is specified below, and define
\begin{equation*}
 A_j=D_h+iS_j,\qquad B_j=D_h-iS_j.
\end{equation*}
Their principal symbols are $\mathfrak a_j$ and $\mathfrak b_j$,
respectively. Microlocally,
\begin{align}\label{car:eq:factorization}
 Q_{j,h}&=A_jB_j+h\kappa R_{j,1}+h^2R_{j,0},
 &&p_j>f/4,\notag\\
 P_h&=A_1A_2B_1B_2+h\kappa R_3+h^2R_2,
 &&p_1,p_2>f/4,\notag\\
 [B_1,B_2]&=h\kappa\widetilde R_1+h^2\widetilde R_0,
 &&p_1,p_2>f/4.
\end{align}
The last subscript denotes the total semiclassical order in
\eqref{car:eq:word-class}; $j$ labels the factor $Q_{j,h}$.
All remainder seminorms are polynomially bounded in $\eta^{-1}$,
and localization errors satisfy \eqref{car:eq:exterior}.

The phase derivatives cancel in
\begin{equation}\label{car:eq:exact-normal-commutator}
 [D_h,A_h]=e^{\phi/h}h^3[D_n,A]e^{-\phi/h}.
\end{equation}
With $\zeta=\xi'+i\nabla'\phi$ and $a=g^{kl}\zeta_k\zeta_l$,
the principal root commutator is
\begin{equation}\label{car:eq:root-commutator-symbol}
 -ih\left(\partial_ns_j-i\nabla'\phi_n\cdot
                              \partial_{\xi'}s_j\right)
 =-ih\frac{(\partial_ng^{kl})\zeta_k\zeta_l}{2s_j}.
\end{equation}
This belongs to $h\kappa S^1_{1,0}$.

The tangential Weyl symbol of $A_h$ is $a+ha_1$ modulo
$h^2S^0_{1,0}$. Since $\{s_j,s_j\}=0$, taking
$s_{j,0}=a_1/(2s_j)$ gives
\begin{equation*}
 S_j^2=A_h+(-1)^j\theta^2+h^2R_{j,0},\qquad
 [D_h,S_j]=h\kappa R_{j,1}+h^2\widetilde R_{j,0}.
\end{equation*}
Since $s_1,s_2$ are functions of the same scalar $a$,
\begin{equation*}
 \{s_1,s_2\}=0,\qquad [S_1,S_2]=h^2R_0.
\end{equation*}
Expansion of $A_jB_j$ and commutation of $B_1$ with $A_2$ give
\eqref{car:eq:factorization}. The denominators involve only $s_j$,
so the bounds are uniform as $\theta\to0$.

\par\medskip\noindent\textbf{Step 3. Boundary energy estimates.}\par

In each of the four intersections, choose finitely many real cutoffs
$\chi_0,\chi_1,\ldots$ with the seminorm bounds of
\eqref{car:eq:exterior} and
\begin{equation*}
 0\le\chi_j\le1,\qquad
 \operatorname{dist}\bigl(\operatorname{supp}\chi_j,
             \operatorname{supp}(1-\chi_{j+1})\bigr)\ge c\eta^M
 \quad(x_n\ge0).
\end{equation*}
Set
\begin{equation*}
 v_0=\operatorname{Op}_h(\chi_0)v,\qquad
 w=\mathcal D v_0,\qquad \operatorname{ord}\mathcal D\le r,
\end{equation*}
with $\mathcal D$ as in \eqref{car:eq:word-class}. By
\eqref{car:eq:exterior}, localization errors throughout this step
are $O(h^6\eta^{-K}\|v\|_{H_h^4(X)})$.

Let $C=hD_n-\operatorname{Op}_h\rho$ be a first-order factor on a
neighborhood of $\operatorname{supp}\chi_0\subset\{|\xi'|\le2L\}$.
For $s\in\mathbb N_0$, $r+s+1\le4$, and $h\le c\eta^M$,
$\operatorname{Im}\rho\le-c\eta^b$ implies
\begin{equation}\label{car:eq:scalar-perfect}
 \|w\|_{H_h^{s+1}}
 \le C\eta^{-K}\bigl(\|Cw\|_{H_h^s}+h^6\|v\|_{H_h^4(X)}\bigr),
\end{equation}
without a boundary condition. If instead
\begin{equation*}
 \begin{gathered}
  H=\xi_n-\operatorname{Re}\rho,\qquad K_0=-\operatorname{Im}\rho,\\
  \alpha K_0^2+\{H,K_0\}\ge c\eta^b,\qquad
  0\le\alpha\le C\eta^{-b'},\qquad w(0)=0,
 \end{gathered}
\end{equation*}
with $\alpha$ constant, then
\begin{equation}\label{car:eq:scalar-half}
 \sqrt h\,\|w\|_{H_h^{s+1}}
 \le C\eta^{-K}\bigl(\|Cw\|_{H_h^s}+h^6\|v\|_{H_h^4(X)}\bigr).
\end{equation}

We use the energy argument of
\cite[Lemmas~B.1 and~B.6, Proposition~B.5]{LRR}, keeping track of the
$\eta$-dependence by localization. Extend $\rho$ smoothly from a
neighborhood of $\operatorname{supp}\chi_2$, using a further nested
cutoff; its required seminorms remain polynomial in $\eta^{-1}$.
The positive symbols used below are
$a=-\operatorname{Im}\rho$ in the first case and
$a=\alpha K_0^2+\{H,K_0\}$ in the second. The latter is independent
of $\xi_n$, so both are tangential symbols. On the support of
$\chi_2$ they satisfy $a\ge c\eta^b$. Set
\begin{equation*}
 a_{\rm ext}=\chi_2a+(1-\chi_2)c\eta^b\ge c\eta^b.
\end{equation*}
This agrees with $a$ near $\operatorname{supp}\chi_1$ and preserves
the polynomial seminorm bounds.
By \eqref{car:eq:exterior}, replacing $a$ by $a_{\rm ext}$ in a
quadratic form on $w=\mathcal Dv_0$ costs at most
$C_Nh^N\eta^{-K_N}\|v\|_{H_h^4(X)}\|w\|_2$ for every fixed $N$.
Sharp G\aa rding therefore gives the local lower bound
$c\eta^b-Ch\eta^{-K}$, up to this localization error. Only finitely
many seminorms enter the estimates below; their exponents are fixed
before choosing $M$ in $h\le c\eta^M$.
With $(u,v)$ linear in $u$, integration by parts gives
\begin{equation}\label{car:eq:negative-root-energy}
 \operatorname{Im}(Cw,w)
 =\frac h2\|w(0)\|_2^2
  +(\operatorname{Op}_h(-\operatorname{Im}\rho)w,w)
  +O(h\eta^{-K})\|w\|_2^2.
\end{equation}
For $w(0)=0$ and $h\alpha\le1/2$, the positive commutator and sharp
G\aa rding give
\begin{equation}\label{car:eq:positive-commutator-energy}
 \|Cw\|_2^2\ge
 (ch\eta^b-Ch^2\eta^{-K})\|w\|_2^2
        -C_Nh^{N+1}\eta^{-K_N}\|v\|_{H_h^4(X)}\|w\|_2.
\end{equation}
Choose $N\ge6$ and $c_0>0$ sufficiently small. Then
\begin{equation*}
 \left.
 \begin{gathered}
  M>\max\{K+b,b'\},\\
  h\le c_0\eta^M
 \end{gathered}
 \right\}
 \quad\Longrightarrow\quad
 Ch\eta^{-K}\le\tfrac c2\eta^b,\qquad h\alpha\le\tfrac12.
\end{equation*}
G\aa rding and Young's inequality yield the $L^2$ bounds; enlarging
$M$ absorbs the order-$h$ subprincipal terms. For
$\Lambda=\langle hD'\rangle$,
\begin{equation*}
 \begin{aligned}
  \|[C,\Lambda^s]w\|_2
   &\le Ch\eta^{-K}\|w\|_{H_h^s(\mathrm{tan})},\\
  w(0)=0&\ \Longrightarrow\ (\Lambda^s w)(0)=0,\\
  hD_nw&=Cw+\operatorname{Op}_h(\rho)w.
 \end{aligned}
\end{equation*}
Absorption and induction give
\eqref{car:eq:scalar-perfect}--\eqref{car:eq:scalar-half}.

On $\operatorname{supp}\chi_0\subset\{p_j<f/2\}$,
\begin{equation*}
 \operatorname{Im}(-if\pm is_j)\le-f/2\le-c\eta,
\end{equation*}
and, for $h\le c\eta^M$, without a boundary condition,
\begin{equation}\label{car:eq:block-energy}
 \|\mathcal D v_0\|_{H_h^{s+2}}
 \le C\eta^{-K}\bigl(\|Q_{j,h}\mathcal D v_0\|_{H_h^s}+h^6\|v\|_{H_h^4(X)}\bigr),
 \qquad r+s+2\le4.
\end{equation}

Indeed, set
\begin{equation*}
 H_c=\begin{pmatrix}0&1\\ f^2-z_j&-2if\end{pmatrix},
 \qquad K_c=-iH_c.
\end{equation*}
Schur triangularization gives, also at repeated roots,
\begin{equation*}
 \operatorname{spec}K_c\subset\{\operatorname{Re}\lambda\le-c\eta\},
 \qquad \|e^{tK_c}\|\le C(1+t)e^{-c\eta t}\quad(t\ge0).
\end{equation*}
Hence
\begin{equation}\label{car:eq:lyapunov}
 S=\int_0^\infty e^{tK_c^*}e^{tK_c}\,dt,
 \qquad K_c^*S+SK_c=-I,
 \qquad cI\le S\le C\eta^{-3}I.
\end{equation}
Differentiation yields
\begin{equation*}
 \|\partial^\alpha S\|\le C_\alpha\eta^{-3(|\alpha|+1)}.
\end{equation*}

With $\operatorname{supp}\chi_2\Subset\{p_j<f/2\}$, matrix G\aa rding gives
\begin{equation*}
 \begin{gathered}
  S_{\rm ext}=\chi_2^2S+(1-\chi_2^2)I,\\
  \operatorname{Op}_h(S_{\rm ext})\ge(c_0-Ch\eta^{-K})I\ge c'I
  \qquad(x_n\ge0,\ h\le c\eta^M).
 \end{gathered}
\end{equation*}
For $W=(w,hD_nw)^t$ and $H_{c,h}=\operatorname{Op}_h(H_c)$, insert
a cutoff equal to one near $\operatorname{supp}\chi_1$. On its bounded
tangential-frequency support, the subprincipal terms have order zero.
For each fixed $N$, this gives microlocally
\begin{equation*}
 \begin{gathered}
 (hD_n-H_{c,h})W=(0,Q_{j,h}w)^t+hR_hW+F_N,\\
 \operatorname{ord}R_h\le0,\qquad
 \|R_h\|_{L^2\to L^2}\le C\eta^{-K},\\
 \|F_N\|_2\le C_Nh^N\eta^{-K_N}\|v\|_{H_h^4(X)}.
 \end{gathered}
\end{equation*}
By \eqref{car:eq:exterior}, choosing $N$ and then $M$ sufficiently
large absorbs $F_N$ into the $h^6$ remainder below.
Pairing with $\operatorname{Op}_h(S_{\rm ext})W$ and using
$i(SH_c-H_c^*S)=I$ and \eqref{car:eq:exterior} yields
\begin{equation}\label{car:eq:block-boundary}
 h(\operatorname{Op}_h(S_{\rm ext}(0))W(0),W(0))\ge0,
\end{equation}
\begin{equation*}
 (1-Ch\eta^{-K})\|W\|_2^2
 \le C\eta^{-K}\bigl(\|Q_{j,h}w\|_2+h^6\|v\|_{H_h^4(X)}\bigr)\|W\|_2.
\end{equation*}
Using
\begin{equation*}
 \begin{aligned}
  Ch\eta^{-K}&\le\tfrac12,\\
  (hD_n)^2w&=Q_{j,h}w-2ifhD_nw
   -\bigl(A_h+(-1)^j\theta^2+h\partial_nf-f^2\bigr)w,
 \end{aligned}
\end{equation*}
tangential commutation and induction give \eqref{car:eq:block-energy}.

For $|\xi'|\ge L$,
\begin{equation*}
 \operatorname{Re}z_j\ge c|\xi'|^2,
 \qquad |\operatorname{Im}z_j|\le C|\xi'|,
 \qquad p_j\ge c\langle\xi'\rangle,
 \qquad |\ell_j|\le C.
\end{equation*}
The principal normal roots $\rho_A=-if-is_j$ and $\rho_B=-if+is_j$ satisfy
\begin{equation}\label{car:eq:high-root-gaps}
 \operatorname{Im}\rho_A=-(f+p_j)\le-c\langle\xi'\rangle,
 \qquad
 \operatorname{Im}\rho_B=p_j-f\ge c\langle\xi'\rangle.
\end{equation}
Pairing with $\Lambda w$ as in \cite[Lemmas~B.1 and~B.4]{LRR} gives
\eqref{car:eq:scalar-perfect} for $C=A_j$, and for $C=B_j$ if
$w(0)=0$; higher orders follow as above.
\par\medskip\noindent\textbf{Step 4. The enhanced local estimate.}\par\nobreak

First let $v$ be compactly supported in a normalized boundary patch
with $0\le x_n\le c_b'\eta/2$ and $v=D_nv=0$ on $x_n=0$.
Choose all enlarged cutoff supports in $0\le x_n<c_b'\eta$.
We claim that, for fixed $M$ and $c,C>0$,
\begin{equation}\label{car:eq:local-enhanced}
 \mathcal E_h(v):=h\|v\|_{H_h^4}
       +\sqrt h\sum_{j=1}^2\|Q_{j,h}v\|_{H_h^2}
 \le C\eta^{-M}\|P_hv\|_2,
\end{equation}
provided
\begin{equation}\label{car:eq:local-thresholds}
 h\le c\eta^M,\qquad \kappa\le c\eta^M,
 \qquad \theta^2\le c_*\eta.
\end{equation}
Put $v_0=\operatorname{Op}_h(\chi_0)v$.
Tangential localization preserves the clamped boundary conditions:
\begin{equation*}
 v_0(0)=D_nv_0(0)=0.
\end{equation*}
If $p_j>f/4$ on $\operatorname{supp}\chi_0$, then $(B_jv_0)(0)=0$.

First assume that $|\xi'|\le2L$ and $\operatorname{Re}s_j=p_j>f/4$
for $j=1,2$ on $\operatorname{supp}\chi_0$. Define
\begin{equation*}
 E_0=\|P_hv_0\|_2+h\kappa\|v_0\|_{H_h^3}
                         +h^2\|v_0\|_{H_h^2}+h^6\|v\|_{H_h^4(X)}.
\end{equation*}
Apply \eqref{car:eq:scalar-perfect} successively with
$(w,s)=(A_2B_1B_2v_0,0)$ and $(B_1B_2v_0,1)$. By
\eqref{car:eq:factorization},
\begin{equation}\label{car:eq:lower-concatenation}
 \|B_1B_2v_0\|_{H_h^2}\le C\eta^{-K}E_0.
\end{equation}
Applying \eqref{car:eq:scalar-half} with
$(w,s)=(B_2v_0,2)$ and $(v_0,3)$ gives
\begin{equation}\label{car:eq:upper-concatenation}
 \sqrt h\|B_2v_0\|_{H_h^3}
                 +h\|v_0\|_{H_h^4}\le C\eta^{-K}E_0.
\end{equation}
The commutator in \eqref{car:eq:factorization} permits exchanging
$B_1,B_2$. The individual factorizations then give
\begin{equation*}
 \mathcal E_h(v_0)
 \le C\eta^{-K}\bigl(\|P_hv_0\|_2+
                   (\kappa+h)\mathcal E_h(v_0)+h^6\|v\|_{H_h^4(X)}\bigr).
\end{equation*}
Indeed,
\begin{equation*}
 h\kappa\|v_0\|_{H_h^3}+h^2\|v_0\|_{H_h^2}
 \le(\kappa+h)\mathcal E_h(v_0).
\end{equation*}
Increasing $M$ absorbs this term, yielding
\begin{equation}\label{car:eq:piece-enhanced}
 \mathcal E_h(v_0)
 \le C\eta^{-K}\bigl(\|P_hv_0\|_2+h^6\|v\|_{H_h^4(X)}\bigr).
\end{equation}

Next assume that $|\xi'|\le2L$ and $p_b<f/2$ on
$\operatorname{supp}\chi_0$ for some $b\in\{1,2\}$.
Apply \eqref{car:eq:block-energy} to $Q_{k,h}v_0$, where
$\{b,k\}=\{1,2\}$. Exact commutation of the factors yields
\begin{equation*}
 \|Q_{k,h}v_0\|_{H_h^2}
          \le C\eta^{-K}(\|P_hv_0\|_2+h^6\|v\|_{H_h^4(X)}).
\end{equation*}
For the other factor, apply \eqref{car:eq:block-energy} if $p_k<f/2$
on $\operatorname{supp}\chi_0$, and combine
\eqref{car:eq:scalar-perfect} and \eqref{car:eq:scalar-half} if
$p_k>f/4$ there. In either case,
\begin{equation*}
 \sqrt h\|v_0\|_{H_h^4}
          \le C\eta^{-K}(\|Q_{k,h}v_0\|_{H_h^2}+h^6\|v\|_{H_h^4(X)}).
\end{equation*}
Absorb the factorization errors by \eqref{car:eq:local-thresholds}.
Since $Q_{b,h}-Q_{k,h}=\pm2\theta^2$, these bounds also
control $\sqrt h\|Q_{b,h}v_0\|_{H_h^2}$, proving
\eqref{car:eq:piece-enhanced} on all three intersections containing
a stable factor.
If $\operatorname{supp}\chi_0\subset\{|\xi'|\ge L\}$, the estimates
following \eqref{car:eq:high-root-gaps} give
\eqref{car:eq:piece-enhanced}.

Choose the partition symbols to satisfy
\begin{equation}\label{car:eq:linear-partition}
 \sum_{a=1}^{J_\eta}\chi_a(x,\xi')=\sigma(x),
 \qquad \sigma=1\text{ near }\operatorname{supp}v,
 \qquad J_\eta\le C\eta^{-K}.
\end{equation}
Choose one member supported in $\{|\xi'|\ge L\}$ and the others
in $\{|\xi'|\le2L\}$, subordinate to the four intersections above,
with nested cutoff symbols as in Step~3.
Linearity of the common quantization gives
\begin{equation}\label{car:eq:exact-reconstruction}
 v=\sum_av_a,
 \qquad Q_{j,h}v=\sum_aQ_{j,h}v_a,
 \qquad v_a=\operatorname{Op}_h(\chi_a)v,
\end{equation}
since $\operatorname{Op}_h(\sigma)$ is multiplication by $\sigma$.
Hence $\mathcal E_h(v)\le\sum_a\mathcal E_h(v_a)$.

Let $X$ be multiplication by a spatial cutoff function or
$X=\operatorname{Op}_h(\chi)$ for a phase-space cutoff symbol
$\chi\in S^0_{1,0}$. The identity
\begin{equation}\label{car:eq:exact-factor-commutator}
 [P_h,X]=[Q_{1,h},X]Q_{2,h}+[Q_{2,h},X]Q_{1,h}
                          +[Q_{1,h},[Q_{2,h},X]]
\end{equation}
and the orders $[Q_{j,h},X]\in h\Psi_h^1$,
$[Q_{1,h},[Q_{2,h},X]]\in h^2\Psi_h^2$ give
\begin{align}\label{car:eq:commutator-margin}
 \|[P_h,X]v\|_2
 &\le C\eta^{-K}\left(
          h\sum_j\|Q_{j,h}v\|_{H_h^1}+h^2\|v\|_{H_h^2}\right)
          \notag\\
 &\le C\eta^{-K}(\sqrt h+h)\mathcal E_h(v).
\end{align}
Using $P_hv_a=\operatorname{Op}_h(\chi_a)P_hv+
[P_h,\operatorname{Op}_h(\chi_a)]v$ and
$h^6\|v\|_{H_h^4(X)}\le h^5\mathcal E_h(v)$, sum \eqref{car:eq:piece-enhanced}
over the partition. Polynomial multiplicity and operator bounds give
fixed $a_0,b_0$ such that
\begin{equation}\label{car:eq:local-final-absorption}
 \mathcal E_h(v)\le C\eta^{-a_0}\|P_hv\|_2
      +C\eta^{-b_0}(\sqrt h+h+h^5)\mathcal E_h(v).
\end{equation}
Choose $M>2b_0+2$ satisfying all previous restrictions. For sufficiently
small $c$ in \eqref{car:eq:local-thresholds},
\begin{equation*}
 C\eta^{-b_0}(\sqrt h+h+h^5)\le\tfrac12,
\end{equation*}
which proves \eqref{car:eq:local-enhanced} for functions supported
in this collar.

\par\medskip\noindent\textbf{Step 5. Interior patches.}\par

In Euclidean coordinates, the real and imaginary Weyl symbols of
$Q_{j,h}$ are
\begin{equation*}
 a=|\xi|^2-|\nabla\phi|^2+(-1)^j\theta^2,\qquad
 b=2\xi\cdot\nabla\phi.
\end{equation*}
By \eqref{car:eq:second-completion} and ellipticity at large frequencies,
\begin{equation*}
 C\eta^{-1}a^2+\{a,b\}\ge c\eta\langle\xi\rangle^4.
\end{equation*}
For $w$ supported in an interior patch and $h\le c\eta$, symbolic
composition and sharp G\aa rding \cite[Theorems~9.5 and~9.11]{Zworski}
applied to the positive-commutator identity give
\begin{equation*}
 \|Q_{j,h}w\|_2^2\ge
 (ch\eta-Ch^2\eta^{-K})\|w\|_{H_h^2}^2.
\end{equation*}
Commuting $(hD)^\beta$, $|\beta|\le s$, and absorbing the
order-$h$ errors for $h\le c\eta^M$ yield
\begin{equation}\label{car:eq:interior-second}
 \sqrt h\|w\|_{H_h^{s+2}}
       \le C\eta^{-K}\|Q_{j,h}w\|_{H_h^s},\qquad s=0,1,2.
\end{equation}
Apply \eqref{car:eq:interior-second} with
$(j,w,s)=(1,Q_{2,h}v,0)$ and $(2,v,2)$, then exchange the factors
to obtain \eqref{car:eq:local-enhanced} for interior support.

To treat a general $v$ compactly supported in a normalized boundary
patch, choose a fixed smooth function $\chi:[0,\infty)\to[0,1]$
equal to one on $[0,1/4]$ and zero on $[1/2,\infty)$. Set
\begin{equation*}
 \chi_b(x_n)=\chi\!\left(\frac{x_n}{c_b'\eta}\right),\qquad
 v_b=\chi_bv,\qquad v_i=(1-\chi_b)v.
\end{equation*}
Then $v_b$ has the support required in Step~4, while $v_i$ vanishes
near the boundary and can be extended by zero for the interior
estimate. Both pieces retain the clamped conditions, and
$|\partial^\alpha\chi_b|\le C_\alpha\eta^{-|\alpha|}$.
Using $v=v_b+v_i$, applying the two estimates, and invoking
\eqref{car:eq:commutator-margin} for these spatial cutoffs give
\begin{equation*}
 \mathcal E_h(v)
 \le C\eta^{-K}\|P_hv\|_2
     +C\eta^{-K'}(\sqrt h+h)\mathcal E_h(v).
\end{equation*}
Using the uniform coordinate and gauge bounds, increase $M>2K'+2$
and decrease $c$ in $h\le c\eta^M$ to absorb the last term. This proves
\eqref{car:eq:local-enhanced} on the whole normalized patch.
\par\medskip\noindent\textbf{Step 6. Rescaling and summation.}\par
\nobreak\smallskip\noindent
With the expansion and derivative orders fixed, denote the exponent
in \eqref{car:eq:local-enhanced} by $M_0$. Choose $R_0$ so that
\begin{equation}\label{car:eq:uniform-curvature-choice}
 CR_0\sup_{0<\rho\le4}
       \rho\bigl(\log(r_c/\rho)\bigr)^{2M_0}\le c.
\end{equation}
Since $\rho|\log\rho|^{2M_0}\to0$, this choice gives
$\kappa\le c\eta^{M_0}$ and, by Lemma~\ref{car:lem:phase},
$\partial_n\phi\ge c_b\eta$ for $0\le x_n<c_b'\eta$ on each
normalized annulus. Decrease $R_0$ to keep the boundary caps in the
fixed tubular neighborhood.

If $\operatorname{dist}(x_0,\partial\Omega)>4R$, use the radial
phase and Step~5. Otherwise $B_{4R}(x_0)$ lies within distance $8R$
of the boundary, so the boundary coordinates apply uniformly.

On an annulus $r\sim\rho$, put $y=\rho z$, $h=\tau^{-1}$,
$\mu=k\rho$, and $v=e^{\tau\Phi}z_0$, where $z_0$ denotes
the test function. The exact scaling identities are
\begin{equation}\label{car:eq:scaling-identities}
 Q_{j,h,z}=h^2\rho^2Q_{j,\tau,y},\qquad
 P_{h,z}=h^4\rho^4P_{\tau,y}.
\end{equation}
Multiply the local enhanced estimate by $h^{-4}\rho^{n/2}$.
Its left side is equivalent to the sum of
\begin{equation}\label{car:eq:scaled-norms}
 \left(\sum_{j=0}^4\tau^{6-2j}\rho^{2j}
                         \|D_y^jv\|_2^2\right)^{1/2},
 \qquad
 \left(\sum_{j=0}^2\tau^{3-2j}\rho^{4+2j}
                    \|D_y^jQ_{\epsilon,\tau}v\|_2^2\right)^{1/2},
 \quad \epsilon=\pm.
\end{equation}
The right side becomes $C\eta^{-M_0}\|\rho^4P_\tau v\|_2$.
If $F$ is the boundary coordinate map and
$\widetilde v=m_0^{-1}F^*v$, then
\begin{equation*}
 \widetilde Q_j\widetilde v=m_0^{-1}F^*(Q_jv),\qquad
 \widetilde P\widetilde v=m_0^{-1}F^*(Pv).
\end{equation*}
The fixed-order bounds on $F^{\pm1}$ and $m_0^{\pm1}$ make the norms in
\eqref{car:eq:scaled-norms} uniformly equivalent in both coordinates.

Choose a dyadic partition $\sum_\rho\chi_\rho=1$ near
$\operatorname{supp}v$, with bounded overlap and
\begin{equation*}
 \operatorname{supp}\chi_\rho\subset\{c\rho<r<C\rho\}\cap\{r<4\},
 \qquad |D^\alpha\chi_\rho|\le C_\alpha\rho^{-|\alpha|}.
\end{equation*}
Subordinate the boundary and interior partitions from Steps~4--5.
Their enlarged supports retain bounded dyadic overlap, and their
multiplicities and seminorms are polynomial in
$t\le T_*+C$. The powers of $\rho$ in
\eqref{car:eq:scaled-norms} cancel the cutoff derivatives.

For a dyadic spatial cutoff, \eqref{car:eq:exact-factor-commutator}
gives in the original scaled variables
\begin{align*}
 \|r^4[P_\tau,\chi]v\|_2\le C\bigg(&
 \sum_\epsilon\|r^3DQ_{\epsilon,\tau}v\|_2
 +\tau\sum_\epsilon\|r^2Q_{\epsilon,\tau}v\|_2\\
 &+\|r^2D^2v\|_2+\tau\|rDv\|_2+\tau^2\|v\|_2\bigg)
 \le C\tau^{-1/2}\mathcal E_\tau(v).
\end{align*}
The remaining cutoff derivatives contribute fixed powers of $T_*$.
The commutators used to reconstruct $Q_{j,\tau}v$ satisfy the same
bound. Reconstruction and square summation give
\begin{equation}\label{car:eq:annular-absorption}
 \mathcal E_\tau(v)
 \le CT_*^{A_0}\|r^4P_\tau v\|_2
       +CT_*^{B_0}\tau^{-1/2}\mathcal E_\tau(v)
\end{equation}
for fixed $A_0,B_0$.

All local restrictions hold for
\begin{equation*}
 \tau\ge C(1+k+T_*^{K_0}).
\end{equation*}
Take $M\ge\max\{A_0,2B_0,K_0\}$ and enlarge $C$. Then
\begin{equation*}
 \tau\ge C(1+k+T_*^M)
 \quad\Longrightarrow\quad CT_*^{B_0}\tau^{-1/2}\le\tfrac12.
\end{equation*}
Absorption in \eqref{car:eq:annular-absorption} proves
\eqref{car:eq:annular-enhanced}. Since
\begin{equation*}
 \tau^3\|v\|_{L^2(\Omega_{x_0,R})}\le\mathcal E_\tau(v),
 \qquad P_\tau v=e^{\tau\Phi}(\Delta^2-k^4)w,
\end{equation*}
we obtain \eqref{car:eq:annular-scalar}.
\end{proof}

\section{Boundary doubling}\label{sec4}

Fix $L\ge8$. For $x\in\overline\Omega$ and $r>0$, define
\begin{equation}\label{eq:doubling-definition}
 H_x(r)=\int_{\Omega\cap B_r(x)}u(y)^2\,dy,
 \qquad
 D_L(x,r)=\log\frac{H_x(Lr)}{H_x(r)}.
\end{equation}
\begin{theorem}[Local doubling]\label{dbl:thm:linear}
Fix $L\geq8$. For every nonzero eigenfunction $u$ of \eqref{main:problem},
there exist $R_H,c_H,C_H>0$, depending only on $\Omega,n,L$, such that
\begin{equation}\label{dbl:linear}
 D_L(x,r)\le C_H\bigl(1+D_L(x,R)+kR\bigr),
 \qquad x\in\overline\Omega,\quad
 0<r\le c_HR,\quad R\le R_H.
\end{equation}
\end{theorem}

\begin{corollary}[Global mass and doubling]\label{tree:global-doubling}
      There is a fixed $\rho>0$ such that
\begin{equation}\label{tree:root-seed}
 H_x(\rho)\ge e^{-C_\Omega(1+k)}\|u\|_{L^2(\Omega,dy)}^2,
 \qquad x\in\overline\Omega,
\end{equation}
and
 \begin{equation}\label{tree:finite-index}
 \sup_{x\in\overline\Omega,\ r>0}D_L(x,r)\le C_\Omega(1+k).
 \end{equation}
The constants may depend on the fixed ratio $L$.
\end{corollary}

Write $\mathbb R^n_+=\{y_n>0\}$, let $x=(0,d)$ with $d\ge0$, and put
\begin{equation*}
 r=|y-x|,\qquad t=\log r,\qquad \omega=(y-x)/r,\qquad
 d\gamma_x=r^{-n}\,dy=dt\,d\omega.
\end{equation*}
The full cylinder is $\mathcal C=\mathbb R\times\mathbb S^{n-1}$.

\begin{lemma}[Flat clamped singular estimate]\label{dbl:lem:flat}
For $\sigma\ge4$ in $\mathbb Z+\tfrac12$ and smooth $u$ compactly
supported in $\overline{\mathbb R^n_+}\setminus\{x\}$, with
$u=\partial_nu=0$ on $\{y_n=0\}$, one has
\begin{align}
 \sigma\|r^{-\sigma}u\|_{L^2(d\gamma_x)}
 &\le C_n\|r^{4-\sigma}\Delta^2u\|_{L^2(d\gamma_x)},\label{dbl:flat-l2}\\
 \sum_{j=0}^4\sigma^{1-j}\|r^{j-\sigma}D_y^ju\|_{L^2(d\gamma_x)}
 &\le C_n\|r^{4-\sigma}\Delta^2u\|_{L^2(d\gamma_x)}.\label{dbl:flat-h4}
\end{align}
The constants are uniform for $d\ge0$ and independent of the support
radii. The same estimates hold on $\mathbb R^n$ for smooth functions
compactly supported away from the pole.
\end{lemma}

\begin{proof}
For $s>0$, use Duffin's clamped reflection \cite{Duffin}
(see also \cite[Eq.~(3.15)]{MRV}):
\begin{equation}\label{dbl:reflection}
 U(y',s)=u(y',s),\qquad
 U(y',-s)=-u(y',s)+2s\partial_su(y',s)-s^2\Delta u(y',s).
\end{equation}
Set $a_j(y')=\partial_s^j u(y',0)$ for $j=2,3$. The clamped conditions
and \eqref{dbl:reflection} give the Taylor expansions
\begin{equation*}
 \begin{aligned}
 u(y',s)&=\tfrac12a_2(y')s^2+\tfrac16a_3(y')s^3+O(s^4),\\
 U(y',-s)&=\tfrac12a_2(y')s^2-\tfrac16a_3(y')s^3+O(s^4),
 \end{aligned}
\end{equation*}
as $s\downarrow0$. The derivatives through order three match at
$y_n=0$, so $U\in H^4_{\rm loc}$. Reflection increases the distance to $x$,
so $U$ has compact support away from $x$.

Put $f=\Delta^2u$, and denote its positive and reflected negative
copies, extended by zero, by $f_+$ and $f_-$. The identities
$[\Delta,s\partial_s]=2\partial_s^2$ and $[\Delta,s^2]=2+4s\partial_s$
give on the negative side
\begin{equation*}
 \Delta^2U(y',-s)=-f-2\partial_s(sf)-\Delta(s^2f).
\end{equation*}
Since $y_nf_-$, $y_n^2f_-$ and $\partial_{y_n}(y_n^2f_-)$
vanish on $y_n=0$,
\begin{equation}\label{dbl:forcing-distribution}
 \Delta^2U=f_+-f_--2\partial_{y_n}(y_nf_-)-\Delta(y_n^2f_-).
\end{equation}
Set
\begin{equation*}
 \begin{gathered}
 L_a=(\partial_t+a)^2+(n-2)(\partial_t+a)+\Delta_{\mathbb S^{n-1}},
 \qquad P_\sigma=L_{\sigma-2}L_\sigma,\\
 W=r^{-\sigma}U,\qquad G_\pm=r^{4-\sigma}f_\pm,\qquad
 c_d=y_n/r,\qquad D_\sigma=r\partial_{y_n}+(\sigma-3)\omega_n.
 \end{gathered}
\end{equation*}
Equation \eqref{dbl:forcing-distribution} becomes
\begin{equation}\label{dbl:cylinder-forcing}
 P_\sigma W=G_+-G_--2D_\sigma(c_dG_-)-L_{\sigma-2}(c_d^2G_-).
\end{equation}
For $s>0$, set $r_\pm=|(y',\pm s)-x|$. Since $r_-\ge r_+$,
\begin{align}\label{dbl:reflected-weight}
 \|G_-\|_{L^2(d\gamma_x)}^2
 &=\int_{s>0}r_-^{8-2\sigma-n}|f(y',s)|^2\,dy'\,ds\notag\\
 &\le\int_{s>0}r_+^{8-2\sigma-n}|f(y',s)|^2\,dy'\,ds
 =\|G_+\|_{L^2(d\gamma_x)}^2.
\end{align}

On a spherical harmonic of degree $\ell$ and at Fourier frequency $\xi$
in $t$, the multiplier of $P_\sigma$ is
\begin{equation}\label{dbl:multiplier}
 (z-\ell)(z-\ell-2)(z+\ell+n-2)(z+\ell+n-4),
 \qquad z=\sigma+i\xi .
\end{equation}
Since $\sigma\in\mathbb Z+\tfrac12$ and $\sigma\ge4$,
\begin{align*}
 |z-\ell|,\ |z-\ell-2|&\ge\tfrac12,\\
 |(z+\ell+n-2)(z+\ell+n-4)|
 &\ge c_n\bigl(\sigma^2+\xi^2+\ell(\ell+n-2)\bigr),\\
 |(z-\ell)(z+\ell+n-2)|&\ge c_n\sigma.
\end{align*}
Hence
\begin{align*}
 \|P_\sigma^{-1}\|_{L^2(\mathcal C,dt\,d\omega)\to L^2(\mathcal C,dt\,d\omega)}
 &\le C_n\sigma^{-2},\\
 \|P_\sigma^{-1}\|_{H^{-1}_\sigma(\mathcal C)\to L^2(\mathcal C,dt\,d\omega)}
 &\le C_n\sigma^{-1},\\
 \|L_\sigma^{-1}\|_{L^2(\mathcal C,dt\,d\omega)\to L^2(\mathcal C,dt\,d\omega)}
 &\le C_n\sigma^{-1},
\end{align*}
where
\begin{equation*}
 \|F\|_{H^{-1}_\sigma(\mathcal C)}=
 \|(\sigma^2-\partial_t^2-\Delta_{\mathbb S^{n-1}})^{-1/2}F\|_{L^2(\mathcal C,dt\,d\omega)}.
\end{equation*}
On $\mathcal C$,
\begin{equation*}
 r\partial_{y_n}
 =\omega_n\partial_t+(e_n-\omega_n\omega)\cdot\nabla_{\mathbb S^{n-1}}.
\end{equation*}
The sphere coefficients are smooth, so duality gives
\begin{equation*}
 \|D_\sigma F\|_{H^{-1}_\sigma(\mathcal C)}
 \le C_n\|F\|_{L^2(\mathcal C,dt\,d\omega)}.
\end{equation*}
Since $P_\sigma=L_{\sigma-2}L_\sigma$, \eqref{dbl:cylinder-forcing} gives
\begin{equation*}
 W=P_\sigma^{-1}(G_+-G_-)-2P_\sigma^{-1}D_\sigma(c_dG_-)
      -L_\sigma^{-1}(c_d^2G_-).
\end{equation*}
The multiplier bounds and \eqref{dbl:reflected-weight} give
\begin{equation*}
 \|W\|_{L^2(\mathcal C,dt\,d\omega)}
 \le C_n\sigma^{-1}\|G_+\|_{L^2(\mathcal C,dt\,d\omega)}.
\end{equation*}
Restriction to $\mathbb R^n_+$ proves \eqref{dbl:flat-l2}.
On $\mathbb R^n$, use $P_\sigma W=r^{4-\sigma}\Delta^2u$.

Cover each annulus $2^m<r<2^{m+1}$ by balls of radius
$s_m=2^m/(64\sigma)$ centered in the half-space, whose eightfold
enlargements have bounded overlap. On these enlargements,
$r\sim2^m$ and $r^{-\sigma-n/2}$ varies by a bounded factor.
Apply the scaled interior $H^4$ estimate when the twice-enlarged
ball misses the plane; otherwise use boundary half-balls of radii
$3s_m$ and $6s_m$ within the eightfold enlargement. In either case,
\begin{equation*}
 s_m^j\|D_y^ju\|_{L^2(dy)}
 \le C\bigl(\|u\|_{L^2(dy)}+s_m^4\|\Delta^2u\|_{L^2(dy)}\bigr),
 \qquad 0\le j\le4,
\end{equation*}
with the right-hand norms taken over the larger set. Multiplying by
the local weight and summing the squared estimates gives
\begin{equation}\label{dbl:elliptic-parameter}
 \sum_{j=0}^4\sigma^{-j}\|r^{j-\sigma}D_y^ju\|_{L^2(d\gamma_x)}
 \le C_n\bigl(\|r^{-\sigma}u\|_{L^2(d\gamma_x)}+
               \sigma^{-4}\|r^{4-\sigma}\Delta^2u\|_{L^2(d\gamma_x)}\bigr).
\end{equation}
Combining this with \eqref{dbl:flat-l2} proves \eqref{dbl:flat-h4}.
The whole-space argument uses interior estimates on all balls.
\end{proof}

Let
\begin{equation*}
 \mathcal P=\Delta_y^2+\sum_{|\alpha|\le4}e_\alpha(y)\partial_y^\alpha
\end{equation*}
be uniformly elliptic on a flattened cap or an interior patch. For
$r=|y-x|\le a\le1$ suppose
\begin{equation}\label{dbl:coefficient-assumptions}
 |e_\alpha(y)|\le Kr\quad(|\alpha|=4),\qquad
 |e_\alpha(y)|\le K\quad(|\alpha|\le3).
\end{equation}
Assume uniform coefficient bounds and local elliptic estimates,
clamped or interior as appropriate. For a smooth clamped function $u$
supported in $0<r<a$, set
\begin{equation*}
 S=\sum_{j=0}^4\sigma^{1-j}\|r^{j-\sigma}D_y^ju\|_{L^2(d\gamma_x)}.
\end{equation*}
By \eqref{dbl:coefficient-assumptions},
\begin{align}\label{dbl:perturbation-error}
 \|r^{4-\sigma}(\mathcal P-\Delta^2)u\|_{L^2(d\gamma_x)}
 &\le CKa\sigma^3S,\notag\\
 \|\mu^4r^{4-\sigma}u\|_{L^2(d\gamma_x)}
 &\le(\mu a)^4\sigma^{-1}S.
\end{align}
Absorption in \eqref{dbl:flat-h4} gives
\begin{equation}\label{dbl:perturbed-estimate}
 \sum_{j=0}^4\sigma^{1-j}\|r^{j-\sigma}D_y^ju\|_{L^2(d\gamma_x)}
 \le C\|r^{4-\sigma}(\mathcal P-\mu^4)u\|_{L^2(d\gamma_x)}
\end{equation}
provided
\begin{equation}\label{dbl:absorption-conditions}
 Ka\sigma^3\le c,\qquad (\mu a)^4/\sigma\le c.
\end{equation}

For the Euclidean bi-Laplacian in Fermi coordinates,
$g^{nn}=1$ and $g^{n\alpha}=0$ for $\alpha<n$.
A constant tangential linear change gives
$g^{\alpha\beta}(x)=\delta_{\alpha\beta}$ at the pole and preserves
$\{y_n=0\}$. Thus the principal error vanishes to first order at $x$.
In a chart rescaled by $R$, \eqref{dbl:coefficient-assumptions} holds
with $K\le C_\Omega R$.

\begin{proposition}[Small-scale propagation from one outer mass ratio]
\label{dbl:prop:small-scale}
Fix $L\ge8$ and $0<s\le1/L$. Let $v$ be a smooth solution of
$(\mathcal P-\mu^4)v=0$ in the open half-space part of a flattened
cap that contains a relative neighborhood of
$B_{Ls}(x)\cap\overline{\mathbb R^n_+}$ in
$\overline{\mathbb R^n_+}$, with clamped boundary conditions on the planar part of
the cap boundary.
Assume \eqref{dbl:coefficient-assumptions} and uniform local clamped
elliptic estimates there. Set
\begin{equation*}
 H(\rho)=\int_{B_\rho(x)\cap\{y_n>0\}}|v|^2\,dy,\qquad
 N=\log\frac{H(Ls)}{H(s)},
\end{equation*}
with $H(\rho)>0$ for the radii considered. There are constants
$c_0,c,C>0$, depending only on $n,L$ and the stated uniform estimates,
such that if $\sigma\ge4$ is a half-integer and
\begin{equation}\label{dbl:small-scale-assumptions}
 \sigma\ge C(1+N),\qquad Ks\sigma^3\le c,\qquad \mu s\le c,
\end{equation}
then
\begin{equation}\label{dbl:small-scale-doubling}
 \log\frac{H(L\rho)}{H(\rho)}\le C\sigma
 \qquad (0<\rho\le c_0s).
\end{equation}
On an interior patch, require a neighborhood of $B_{Ls}(x)$
and use full-ball integrals and interior elliptic estimates.
\end{proposition}

\begin{proof}
Use a cutoff $\chi$ supported in $\rho/4<r<3s$ and equal to one
for $\rho/2\le r\le2s$, with derivatives bounded at the respective
scales. The conditions in \eqref{dbl:small-scale-assumptions}
imply \eqref{dbl:absorption-conditions} at $a=3s$ after fixing $c$
sufficiently small. The $L^2$ term of \eqref{dbl:perturbed-estimate}
and local regularity at the two cutoff scales give
\begin{equation}\label{dbl:small-scale-cutoff}
 \|r^{-\sigma}\chi v\|_{L^2(d\gamma_x)}
 \le C\bigl(4^\sigma\rho^{-\sigma-n/2}H(\rho)^{1/2}
             +2^{-\sigma}s^{-\sigma-n/2}H(Ls)^{1/2}\bigr).
\end{equation}
After multiplication by $r^4$, the commutators involve scaled derivatives
of order at most three, estimated on $B_\rho(x)$ and $B_{4s}(x)$
with uniform constants since $\mu s\le c$.

For $\rho\le s$, the region $\rho<r<s$ gives
\begin{equation*}
 H(s)^{1/2}\le H(\rho)^{1/2}
      +s^{\sigma+n/2}\|r^{-\sigma}\chi v\|_{L^2(d\gamma_x)}.
\end{equation*}
Substituting into \eqref{dbl:small-scale-cutoff} and using
$H(Ls)^{1/2}=e^{N/2}H(s)^{1/2}$ and $C2^{-\sigma}e^{N/2}\le1/2$, we obtain
\begin{equation*}
 \|r^{-\sigma}\chi v\|_{L^2(d\gamma_x)}
 \le C4^\sigma\rho^{-\sigma-n/2}H(\rho)^{1/2}.
\end{equation*}
For $c_0\le1/L$, the weight on $\rho<r<L\rho$ then gives
\begin{equation*}
 (L\rho)^{-\sigma-n/2}(H(L\rho)-H(\rho))^{1/2}
 \le C4^\sigma\rho^{-\sigma-n/2}H(\rho)^{1/2},
\end{equation*}
which implies \eqref{dbl:small-scale-doubling}.
\end{proof}

\begin{proof}[Proof of Theorem~\ref{dbl:thm:linear}]
\par\medskip\noindent\textbf{Step 1. Geometry and rescaling.}
Decrease $R_H$ below the radius $R_0$ in Theorem~\ref{car:thm:annular}
so that any pole within $4R$ of the boundary has a Fermi chart containing
all the neighborhoods required inside $B_{4R}(x)$. Normalize the metric
at the pole as above. In the coordinates scaled by $R$, the pulled-back
operator obeys \eqref{dbl:coefficient-assumptions} with
$K\le C_\Omega R$. At a pole whose $4R$-ball is interior, use Euclidean
coordinates and $K=0$. All local elliptic estimates have uniform constants.

Let $G\ge1$ be a common bi-Lipschitz constant for these normalized
charts, and let $m_-,m_+>0$, with $m_-\le m_+$, bound the associated
volume comparisons. Write
\begin{equation*}
 u_R(z)=u(x+Rz),\qquad
 H(t)=\int_{R^{-1}(\Omega-x)\cap B_t(0)}|u_R(z)|^2\,dz
     =R^{-n}H_x(Rt).
\end{equation*}
Let $v$ be the coordinate pullback of $u_R$ and $x_c$ the coordinate pole.
In a boundary chart, set
\begin{equation*}
 H_c(t)=\int_{B_t(x_c)\cap\mathbb R^n_+}|v(y)|^2\,dy;
\end{equation*}
in an interior chart, integrate over $B_t(x_c)$. For balls contained
in the chart,
\begin{equation}\label{dbl:mass-comparison}
 m_-H(t/G)\le H_c(t)\le m_+H(Gt).
\end{equation}
Set
\begin{equation}\label{dbl:ratios}
 L_p=G^2L,\qquad \Gamma=G^2L_p=G^4L,\qquad
 N=D_L(x,R),\qquad \nu=kR,\qquad B=1+N+\nu.
\end{equation}
Set $T_*=\log(r_c/\delta_*)$ and
$\phi_o^+=\sup_{2<r<3}\Phi$. The phase bounds in Section~\ref{sec3} give
\begin{equation}\label{dbl:phase-gap}
 \inf_{0<r\le1}\Phi-\phi_o^+\ge g>0.
\end{equation}

\par\medskip\noindent\textbf{Step 2. The annular range.}
Let an inner cutoff vanish for $r\le\rho/4$ and equal one for
$r\ge\rho/2$, and let an outer cutoff equal one for $r\le2$ and
vanish for $r\ge11/4$. Their product $\chi$ preserves the clamped
boundary conditions and has support with $\delta_*<r<3$.
Apply \eqref{car:eq:annular-scalar} to $\chi u_R$. Scaled local elliptic
estimates on slightly enlarged transition annuli of the cutoffs give
\begin{equation}\label{dbl:bridge-cutoff}
 \begin{aligned}
 \|e^{\tau\Phi}\chi u_R\|_{L^2(\Omega_{x,R},dz)}
 &\le\Pi\left(e^{\tau\phi_i^+(\rho)}H(\rho)^{1/2}
              +e^{\tau\phi_o^+}H(4)^{1/2}\right),
 \\
 1\le\Pi&\le C T_*^p(1+\tau+\nu)^q,
 \end{aligned}
\end{equation}
where $p,q$ are fixed and $\phi_i^+(\rho)=\sup_{\rho/4<r<\rho/2}\Phi$.
On this annulus, $r^4[\Delta^2,\chi]u_R$
is a bounded-coefficient sum of $\rho^jD^ju_R$, $0\le j\le3$;
local regularity bounds its norm by a polynomial in $1+\nu\rho$
times $H(\rho)^{1/2}$. The outer estimates use neighborhoods in $B_4$,
with $H(4)\le H(L)=e^NH(1)$.

For $\Gamma\rho\le1$, set
$\phi_t^-(\rho)=\inf_{\rho<r<\Gamma\rho}\Phi$. Then
\begin{equation}\label{dbl:phase-oscillation}
 \phi_i^+(\rho)-\phi_t^-(\rho)\le C_\Gamma.
\end{equation}
For a fixed $0<\epsilon\le1$, to be chosen below, set
\begin{equation}\label{dbl:parameters}
 \tau=A(1+N+\nu)=AB,\qquad
 \delta=\epsilon\tau^{-3},\qquad \delta_*=\delta/16.
\end{equation}
Since $\nu,N\le\tau/A$ and $T_*=\log(16r_c/\epsilon)+3\log\tau$,
for each $\epsilon$, all sufficiently large $A$ ensure, for every $B\ge1$,
\begin{equation}\label{dbl:large-parameter-conditions}
 \tau\ge C(1+\nu+T_*^M),\qquad
 \log\Pi\le\min\{1,g/4\}\tau,\qquad
 \Pi e^{N/2-g\tau}\le\tfrac12.
\end{equation}
This follows from $(\log\tau)^a/\tau\to0$ for every fixed $a$.
Also require $\delta\le\Gamma^{-1}$.

Fix $\delta\le\rho\le\Gamma^{-1}$ and write
$S=\|e^{\tau\Phi}\chi u_R\|_{L^2(\Omega_{x,R},dz)}$. The weight on $\rho<r<1$ gives
\begin{equation*}
 H(1)^{1/2}\le H(\rho)^{1/2}
                   +e^{-\tau\inf_{0<r\le1}\Phi}S .
\end{equation*}
Substitute into \eqref{dbl:bridge-cutoff}, using
$H(4)^{1/2}\le e^{N/2}H(1)^{1/2}$. By \eqref{dbl:phase-gap} and
\eqref{dbl:large-parameter-conditions},
\begin{equation*}
 S\le4\Pi e^{\tau\phi_i^+(\rho)}H(\rho)^{1/2}.
\end{equation*}
Since $\chi=1$ on $\rho<r<\Gamma\rho$, \eqref{dbl:phase-oscillation} gives
\begin{equation*}
 \frac{H(\Gamma\rho)}{H(\rho)}
 \le1+16\Pi^2e^{2C_\Gamma\tau}
\end{equation*}
and hence
\begin{equation}\label{dbl:bridge}
 \log\frac{H(\Gamma\rho)}{H(\rho)}\le C_1\tau,
 \qquad \delta\le\rho\le\Gamma^{-1}.
\end{equation}
The bound $\log\Pi\le\tau$ makes $C_1$ independent of
$\epsilon$ and the final choice of $A$.

\par\medskip\noindent\textbf{Step 3.  The singular range.}
Proposition~\ref{dbl:prop:small-scale}, with ratio $L_p$, gives
$c_p,c_a,C_0>0$, $c_p\le1/L_p$, such that
\begin{equation*}
 \begin{gathered}
 N_s=\log\frac{H_c(L_ps)}{H_c(s)},\qquad
 \sigma\in\mathbb Z+\tfrac12,\qquad
 \sigma\ge\max\{4,C_0(1+N_s)\},\\
 Ks\sigma^3\le c_a,\qquad \nu s\le c_a
 \end{gathered}
\end{equation*}
imply
\begin{equation}\label{dbl:micro-implication}
 \log\frac{H_c(L_pt)}{H_c(t)}\le C_0\sigma
 \quad(0<t\le c_ps).
\end{equation}
The support enlargement to $3s$ is included in $c_a$.

Choose, independently of $\epsilon$ and $A$,
\begin{equation}\label{dbl:micro-radius}
 K_*\ge\max\{2G,\,2/(Gc_p)\},\qquad s=K_*\delta.
\end{equation}
Since $s\le K_*A^{-3}$, increasing the lower bound on $A$ ensures
$s\le1/L_p$, $s/G\le\Gamma^{-1}$, and that all balls through radius
$L_ps$ lie in the chart. Thus \eqref{dbl:mass-comparison} applies.
Because $s/G\ge\delta$, that comparison and
\eqref{dbl:bridge} give
\begin{align}
 N_s
 &\le\log(m_+/m_-)
       +\log\frac{H(GL_ps)}{H(s/G)}\notag\\
 &=\log(m_+/m_-)
       +\log\frac{H(\Gamma(s/G))}{H(s/G)}
 \le C_2\tau,                                      \label{dbl:micro-index}
\end{align}
where $C_2$ is independent of $\epsilon$ and $A$.

Take the least half-integer $\sigma$ satisfying
$\sigma\ge\max\{4,C_0(1+N_s)\}$. Then
\begin{equation}\label{dbl:sigma-bound}
 4\le\sigma\le C_\sigma\tau,\qquad
 C_\sigma=5+C_0(1+C_2).
\end{equation}
It follows that
\begin{equation}\label{dbl:micro-absorption}
 Ks\sigma^3\le C_\Omega R_HK_*C_\sigma^3\epsilon,\qquad
 \nu s\le K_*\epsilon\tau^{-2}.
\end{equation}
First choose $\epsilon$ so that
$C_\Omega R_HK_*C_\sigma^3\epsilon\le c_a$. Then choose $A$ so that
$K_*\epsilon A^{-2}\le c_a$ and all preceding annular and geometric
restrictions hold. Both choices depend only on the geometry, $n$, and $L$.
By \eqref{dbl:micro-implication},
\begin{equation*}
 \log\frac{H_c(L_pt)}{H_c(t)}\le C_0C_\sigma\tau,
 \qquad 0<t\le c_ps.
\end{equation*}

\par\medskip\noindent\textbf{Step 4. Joining the ranges.}
The two sides of \eqref{dbl:mass-comparison} give
\begin{equation}\label{dbl:return-to-physical}
 \frac{H(L\rho)}{H(\rho)}
 \le\frac{m_+}{m_-}
       \frac{H_c(GL\rho)}{H_c(\rho/G)}
 =\frac{m_+}{m_-}
       \frac{H_c(L_p(\rho/G))}{H_c(\rho/G)}.
\end{equation}
This gives the $C\tau$ bound for $0<\rho\le Gc_ps$.
Since $Gc_ps\ge2\delta$ and $L\le\Gamma$, \eqref{dbl:bridge} covers
the remaining radii:
\begin{equation*}
 \log\frac{H(L\rho)}{H(\rho)}\le C\tau,
 \qquad 0<\rho\le\Gamma^{-1}.
\end{equation*}
Taking $c_H=(2\Gamma)^{-1}$, scaling back by $R$, and recalling
\eqref{dbl:parameters} proves \eqref{dbl:linear}.
\end{proof}

\begin{proof}[Proof of Corollary~\ref{tree:global-doubling}]
Fix $0<R_*\le R_H$, take $c_H\le1$, and put $\rho=c_HR_*/4$.
Normalize $\|u\|_{L^2(\Omega,dy)}=1$. Cover $\overline\Omega$
by finitely many balls $B_s(z_i)$, $1\le i\le M$, with
$z_i\in\Omega$ and $s=\frac14(1-L^{-1})\rho$.
Join $z_i$ and $z_j$ when their covering balls intersect in $\Omega$.
Connectedness of $\Omega$ makes this graph connected.
The cover gives $H_{z_i}(\rho)\ge M^{-1}$ at some center.

For any $z\in\overline\Omega$, write
$a(z)=-\log H_z(\rho)$. Monotonicity gives
$D_L(z,R_*)\le a(z)$. Applying \eqref{dbl:linear} with
$r=\rho/L$ and $R=R_*$ yields
\begin{equation}\label{tree:mass-step}
 -\log H_z(\rho/L)
 =a(z)+D_L(z,\rho/L)
 \le (1+C_H)a(z)+C_H(1+kR_*).
\end{equation}
For adjacent centers, $|z_i-z_j|<2s$ and
$2s+\rho/L<\rho$. Hence
$B_{\rho/L}(z_i)\subset B_\rho(z_j)$, so that
\begin{equation*}
 a(z_j)\le (1+C_H)a(z_i)+C_H(1+kR_*).
\end{equation*}
Iterating from the center with $H_{z_i}(\rho)\ge M^{-1}$ along paths
of at most $M-1$ edges gives
$a(z_i)\le C(1+k)$ at all covering centers. For arbitrary
$x\in\overline\Omega$, choose $z_i$ with $|x-z_i|<s$.
The inclusion $B_{\rho/L}(z_i)\subset B_\rho(x)$ and
\eqref{tree:mass-step} give the same bound for $a(x)$.
This proves \eqref{tree:root-seed} after restoring the normalization.

In particular, $D_L(x,R_*)\le C(1+k)$ for $x\in\overline\Omega$.
Equation \eqref{dbl:linear} gives \eqref{tree:finite-index} for
$0<r\le c_HR_*$. For $r>c_HR_*$,
\begin{equation*}
 D_L(x,r)\le\log\frac{\|u\|_{L^2(\Omega,dy)}^2}{H_x(\rho)}
 \le C(1+k).
\end{equation*}
\end{proof}

\section{Nodal measure and sharpness}\label{sec5}

Write $Z_u=\{x\in\Omega:u(x)=0\}$ and use
analogous notation for other functions.

\begin{lemma}\label{nod:interior}
There are dimensional constants $L_{\rm int}>1$, $c_{\rm int}\in(0,1)$
and $C>0$ with the following property. Suppose $u\in L^2(\Omega)$
is real, solves $\Delta^2u=k^4u$ in $B_{2r}(z)\subset\Omega$ with
$k\ge0$, and is not identically zero in $B_{2r}(z)$. If
$L\ge L_{\rm int}$ and
\begin{equation}\label{nod:interior-doubling}
 D_L(b,s)\le N\qquad (b\in B_r(z),\ 0<s\le r),
\end{equation}
then
\begin{equation}\label{nod:interior-bound}
 \mathcal H^{n-1}(Z_u\cap B_{c_{\rm int}r}(z))
 \le C(N+kr+1)r^{n-1}.
\end{equation}
\end{lemma}

\begin{proof}
Rescale to $r=1$, $z=0$ and put $\mu=kr$, $M=\|u\|_{L^2(B_2)}$.
The lift $e^{e^{i\pi/4}\mu t}u(x)$ and
\eqref{nod:analytic-regularity} give a holomorphic extension with
\begin{equation}\label{nod:complex-bound}
 \sup_{B_\rho^{\mathbb C^n}(0)}|u|\le Ce^{C\mu}M
\end{equation}
for a dimensional $\rho\in(0,1/4)$. Set $\delta=\rho/4$ and
$v(\zeta)=u(\delta\zeta)$, which is holomorphic in a neighborhood
of $\overline{B_2^{\mathbb C^n}(0)}$.

Take $L_{\rm int}=15/\delta$. For $y\in B_{1/10}(0)\subset\mathbb R^n$
and $L\ge L_{\rm int}$, we have $B_2\subset B_{L\delta/5}(\delta y)$.
Thus \eqref{nod:interior-doubling} gives
\begin{equation}\label{nod:anchor-lower}
 \sup_{B_{1/5}(y)}|v|
 \ge |B_{\delta/5}|^{-1/2}
       \|u\|_{L^2(B_{\delta/5}(\delta y))}
 \ge c e^{-N/2}M.
\end{equation}
Together with \eqref{nod:complex-bound}, this gives
\begin{equation*}
 \sup_{B_2^{\mathbb C^n}(0)}|v|
 \le e^{C(N+\mu+1)}\sup_{B_{1/5}(y)}|v|,
 \qquad y\in B_{1/10}(0).
\end{equation*}
By \cite[Proposition~6.7]{DF88},
\begin{equation*}
 \mathcal H^{n-1}(Z_u\cap B_{\delta/20}(0))
 =\delta^{n-1}\mathcal H^{n-1}(Z_v\cap B_{1/20}(0))
 \le C(N+\mu+1).
\end{equation*}
Rescaling back proves \eqref{nod:interior-bound} with
$c_{\rm int}=\delta/20$.
\end{proof}

\begin{proof}[Proof of Theorem~\ref{main:sharp}]
  \par\medskip\noindent\textbf{Step 1. Boundary cubes and doubling indices.}
Put $d=n-1$. For a closed tangential cube $Q$ of side $\ell(Q)$,
in the graph coordinates of Section~\ref{sec2} set
\begin{equation*}
 W(Q)=F(Q\times(0,\ell(Q)]),\qquad
 E(Q)=F(16Q\times[0,16\ell(Q)]) ,
\end{equation*}
Here $16Q$ is the concentric enlargement. Define
\begin{equation}\label{tree:index}
 N(Q)=\sup_{\substack{x\in E(Q)\\0<r\le\ell(Q)}}D_L(x,r),
 \qquad \mu(Q)=k\ell(Q).
\end{equation}
For an integer $A\ge2$, let $Q_i$ be the $A^d$ children of side
$\ell(Q)/A$. Then
\begin{equation}\label{tree:nesting}
 E(Q_i)\subset E(Q),\qquad N(Q_i)\le N(Q),\qquad
 \mu(Q_i)=\mu(Q)/A .
\end{equation}

Choose boundary charts with $|\nabla g|\le1/20$ and margins for the
enlargements below. In coordinates scaled about the boundary point
$p$ over the center of $Q$ by $\ell(Q)$, the base is
$[-1/2,1/2]^d$. Translation and dilation preserve the slope bound,
while graph derivatives of order $j\ge2$ acquire the factor
$\ell(Q)^{j-1}$. Taking the root side lengths at most one gives
uniform smooth bounds on all the rescaled charts.
Use
Lemma~\ref{lem-cauchy} at a fixed fraction of this unit scale so that
its observation patch $\Gamma$ lies over $(-1/4,1/4)^d$ and its
target $B_\rho(y_0)$ lies in $E(Q)$ with $\rho<1$.
Fix $C_{\rm out}$ so that $B_{C_{\rm out}}(p)$ contains the Cauchy
domain and a unit neighborhood of $E(Q)$ in these coordinates.
Choose $L\ge\max(8,L_{\rm int})$ uniformly in the charts so that
$B_{C_{\rm out}}(p)\subset B_{L\rho}(y_0)$. Fix the constants $C_H,c_H,R_H$ in
\eqref{dbl:linear} for this $L$. Take the root cubes small enough
that $\ell(Q)\le R_H$ and the physical enlargements lie in a
single graph chart.

\par\medskip\noindent\textbf{Step 2. Boundary values and child estimates.}
Set $N=N(Q)$ and $\mu=\mu(Q)$. Use the unit-scale coordinates
above for all masses and boundary values, so $\Delta^2u=\mu^4u$.
With $M=H_p(C_{\rm out})^{1/2}>0$,
the definition of $N$ and Lemma~\ref{lem-cauchy} give
\begin{equation}\label{tree:cauchy-use}
 e^{-N/2}M\le \|u\|_{L^2(B_\rho(y_0))}
 \le Ce^{C\mu}M^{1-\alpha}T^\alpha,
 \qquad
 T=\|\Delta u\|_{L^2(\Gamma)}
    +\|\partial_\nu\Delta u\|_{L^2(\Gamma)} .
\end{equation}
Since $\Gamma$ has bounded area, some $p_*\in\Gamma$ satisfies
\begin{equation}\label{tree:trace-point}
 |\Delta u(p_*)|+|\partial_\nu\Delta u(p_*)|
 \ge c\exp\!\left(-\frac{N}{2\alpha}-C\mu\right)M.
\end{equation}
Let $Q_*$ be a child containing the tangential coordinate of $p_*$.
By the uniform Lipschitz bound for $F$,
\begin{equation}\label{tree:child-diameter}
 |x-p_*|\le C/A\qquad(x\in E(Q_*)).
\end{equation}

Choose an integer $m>4C_H/\alpha$ and fix radii
$R_j=\tfrac12 L^{j-m}$, $0\le j\le m$. For $A$ sufficiently
large, $A^{-1}\le c_H R_0$ and
$B_{R_0/2}(p_*)\subset B_{R_0}(x)$ for all $x\in E(Q_*)$.
Boundary regularity \eqref{nod:finite-regularity} in a fixed
graph neighborhood of $p_*$ inside $B_{R_0/2}(p_*)$ bounds the
left-hand side of \eqref{tree:trace-point} by
$C_m(1+\mu)^{b}H_x(R_0)^{1/2}$. Thus
\begin{equation}\label{tree:child-mass}
 \log\frac{M^2}{H_x(R_0)}
 \le \frac{N}{\alpha}+C_m(1+\mu).
\end{equation}
For $r\le A^{-1}$, use \eqref{dbl:linear} at each $R_j$,
$0\le j<m$. The physical radii are $\ell(Q)R_j$ and the spectral
terms are $\mu R_j$. Summation and $H_x(R_m)\le M^2$ give
\begin{align*}
 D_L(x,r)
 &\le C_H\left(1+\frac1m\log\frac{H_x(R_m)}{H_x(R_0)}
             +\frac{\mu}{m}\sum_{j=0}^{m-1}R_j\right)\\
 &\le \frac{C_H}{\alpha m}N+C_m(1+\mu).
\end{align*}
For $K\ge\max(1,4C_m)$ and $N>K(1+\mu)$, taking the supremum gives
\begin{equation}\label{tree:good-child-inequality}
 N(Q_*)\le\tfrac14N+C_m(1+\mu)\le\tfrac12N.
\end{equation}
The choice of $K$ is
independent of subsequent increases of $A$.

Now suppose $\mu\le1$ and $N\le2K$. Write $U=u\circ F$.
Since $U(y,0)=\partial_tU(y,0)=0$,
\begin{equation*}
 U=t^2V,\qquad V(y,t)=\int_0^1(1-b)\partial_t^2U(y,bt)\,db.
\end{equation*}
The boundary $C^{3,\gamma}$ estimate for $U$ near $\Gamma$ gives
$\|V\|_{C^{1,\gamma}}\le CM$ on a smaller neighborhood, for fixed
$0<\gamma<1$. With $a=1+|\nabla g|^2$,
\begin{equation*}
 W:=(\Delta u)\circ F
 =\Delta_y U-2\nabla g\cdot\nabla_y\partial_t U
   +a\partial_t^2U-(\Delta g)\partial_tU.
\end{equation*}
At $t=0$,
\begin{equation*}
 W=2aV,\qquad
 \partial_tW=6a\partial_tV-4\nabla g\cdot\nabla_yV-2(\Delta g)V,
 \qquad
 (\partial_\nu\Delta u)\circ F
 =a^{-1/2}(\nabla g\cdot\nabla_y W-a\partial_tW).
\end{equation*}
Thus \eqref{tree:trace-point} gives
\begin{equation*}
 |V|+|\nabla_{y,t}V|\ge c_KM
 \quad\text{at the coordinate point of }p_*.
\end{equation*}

The $C^{1,\gamma}$ bound gives a rectangular neighborhood of uniform
size on the side $t\ge0$ where either $|V|\ge c_K'M$ or a coordinate
derivative has fixed sign and modulus at least $c_K'M$. In the second case,
write the coordinates as $(s,z)$ with $s$ in that direction.
For any two zeros $(s,z),(s',z')$ in the rectangle, the mean-value
theorem on its coordinate segments gives
\begin{equation*}
 c_K'M|s-s'|
 \le |V(s',z)-V(s',z')|\le CM|z-z'|.
\end{equation*}
Hence the zeros lie on a Lipschitz graph over the $z$ coordinates.
With $K$ fixed, increase $A$ so that $F^{-1}(W(Q_*))$ lies in this
rectangle, using \eqref{tree:child-diameter}. Each projected coordinate,
including $t$ when $s$ is tangential, ranges over an interval of length
$O(A^{-1})$. The projection has measure $O(A^{-d})$, and the graph
and $F$ have bounded Lipschitz constants. Since $Z_U=Z_V$ for $t>0$,
rescaling gives
\begin{equation}\label{tree:terminal}
 \mathcal H^{n-1}(Z_u\cap W(Q_*))
 \le C_K\bigl(\ell(Q)/A\bigr)^d.
\end{equation}
If $|V|\ge c_K'M$, then $Z_u\cap W(Q_*)=\varnothing$, so the same bound holds.
Fix $K$ and $A$ henceforth; $K,A,C_K$ depend only on the geometry and dimension.

\par\medskip\noindent\textbf{Step 3. Summation over boundary cubes.}
Set
\begin{equation}\label{tree:cap-decomposition}
 V(Q)=F(Q\times[\ell(Q)/A,\ell(Q)]),\qquad
 W(Q)=V(Q)\cup\bigcup_i W(Q_i).
\end{equation}
The band $V(Q)$ is at distance at least $c\ell(Q)/A$ from the
physical boundary. Cover it by a bounded number of target balls
from Lemma~\ref{nod:interior}, with radii comparable to
$\ell(Q)/A$. The doubling centers lie in $E(Q)$ and the radii
are at most $\ell(Q)$. Consequently
\begin{equation}\label{tree:band}
 \mathcal H^{n-1}(Z_u\cap V(Q))
 \le C_A\bigl(N(Q)+\mu(Q)+1\bigr)\ell(Q)^d.
\end{equation}

At each cube, estimate $V(Q)$ by \eqref{tree:band} and continue on
the children in
\begin{equation*}
 \mathcal C(Q)=
 \begin{cases}
  \{Q_i:Q_i\ne Q_*\},
    &N(Q)\le K(1+\mu(Q)),\quad\mu(Q)\le1,\\
  \{Q_i\},&\text{otherwise}.
 \end{cases}
\end{equation*}
In the first case, estimate $W(Q_*)$ by \eqref{tree:terminal} and stop
on $Q_*$.

Set $\mathcal F(Q)=N(Q)+\mu(Q)+1$ and
\begin{equation*}
 \varepsilon=\min\left\{
 \frac1{4A^d},\frac{1-A^{-1}}{2(K+1)},A^{-d}
 \right\}>0.
\end{equation*}
Writing $N=N(Q)$ and $\mu=\mu(Q)$, we obtain from
\eqref{tree:nesting}, \eqref{tree:good-child-inequality}, and the
stopping rule that
\begin{equation*}
 \mathcal F(Q)-A^{-d}\sum_{Q_i\in\mathcal C(Q)}\mathcal F(Q_i)
 \ge
 \begin{cases}
  N/(2A^d), & N>K(1+\mu),\\
  (1-A^{-1})\mu, & N\le K(1+\mu),\quad\mu>1,\\
  A^{-d}\mathcal F(Q), & N\le K(1+\mu),\quad\mu\le1.
 \end{cases}
\end{equation*}
Use $\mathcal F(Q)\le2N$ in the first case and
$\mathcal F(Q)\le2(K+1)\mu$ in the second. Since
$\ell(Q_i)^d=A^{-d}\ell(Q)^d$, the choice of $\varepsilon$ gives
\begin{equation}\label{tree:budget}
 \sum_{Q_i\in\mathcal C(Q)}\mathcal F(Q_i)\ell(Q_i)^d
 \le (1-\varepsilon)\mathcal F(Q)\ell(Q)^d
\end{equation}
in all three cases.

By \eqref{tree:band} and \eqref{tree:terminal}, the nodal measure in
$V(Q)$ and, when $Q_*$ is stopped, in $W(Q_*)$ is at most
$C\mathcal F(Q)\ell(Q)^d$.
Since $A^{-j}\ell(Q_0)\to0$, every point of $W(Q_0)$ belongs at finite
depth to some $V(Q)$ or to $W(Q_*)$ for a stopped child $Q_*$.
The closed base cubes and bands
include the artificial faces. Iterating \eqref{tree:budget} gives
\begin{equation}\label{tree:cap-bound}
 \mathcal H^{n-1}(Z_u\cap W(Q_0))
 \le C\sum_{j=0}^\infty(1-\varepsilon)^j
          \mathcal F(Q_0)\ell(Q_0)^d
 =\frac{C}{\varepsilon}\mathcal F(Q_0)\ell(Q_0)^d.
\end{equation}

Choose finitely many root cubes $Q_0$ so that the regions $W(Q_0)$
cover a collar of $\partial\Omega$ in $\Omega$.
By \eqref{tree:finite-index}, $\mathcal F(Q_0)\le C_\Omega(1+k)$
for each root cube. Applying
\eqref{tree:cap-bound} on the collar and
Lemma~\ref{nod:interior} on a finite ball cover of the remaining
compact interior region gives
\begin{equation*}
 \mathcal H^{n-1}(Z_u)\le C_\Omega(1+k).
\end{equation*}
The clamped energy identity and Poincare's inequality give
$k^4\|u\|_2^2=\|\Delta u\|_2^2\ge c_\Omega\|u\|_2^2$, so
$1+k\le C_\Omega k$.
\end{proof}

\begin{example}[Sharpness on the unit ball]\label{tree:optimality}
On the unit ball $\mathbb B^n$, $n\ge2$, let
$Y_m(\omega)=\operatorname{Re}(\omega_1+i\omega_2)^m$ and
$\sigma_m=m(m+n-2)$, $m\ge1$. Minimizing the clamped Rayleigh
quotient in the invariant angular sector $a(r)Y_m(\omega)$ gives
a real eigenfunction $u_m=a_m(r)Y_m(\omega)$. The angular energy
bound $\|\nabla u_m\|_2^2\ge\sigma_m\|u_m\|_2^2$ and a fixed
radial test function supported in an interior annulus yield
\begin{equation*}
 \sigma_m^2\le\lambda_m\le C(1+\sigma_m)^2,
 \qquad \lambda_m\sim m^4.
\end{equation*}
The nodal set contains the sections of the $m$ distinct hyperplanes
defined by $\operatorname{Re}(x_1+ix_2)^m=0$. Their pairwise
intersections have zero $(n-1)$-dimensional measure, so
\begin{equation*}
 \mathcal H^{n-1}(Z_{u_m})\ge c_n m\ge c_n'\lambda_m^{1/4}.
\end{equation*}
\end{example}

\section*{Acknowledgments}
Liu's research is supported by the National Natural Science Foundation of China (No. 12471198).
Yang's research is supported by the National Natural Science Foundation of China (No. 12431018).

\bigskip
\noindent\textsuperscript{1}School of Mathematics and Statistics,
Nanjing University of Science and Technology,
Nanjing 210094, P.R. China\par
\noindent\textit{Email address:} hrliu@njust.edu.cn

\medskip
\noindent\textsuperscript{2}School of Mathematics and Statistics,
Nanjing University of Science and Technology,
Nanjing 210094, P.R. China\par
\noindent\textit{Email address:} tianlong19850812@163.com

\medskip
\noindent\textsuperscript{3}School of Mathematics,
Nanjing University, Nanjing 210093, P.R. China\par
\noindent\textit{Email address:} xpyang@nju.edu.cn

\end{document}